\documentclass[11pt]{article}

\usepackage{latexsym}
\usepackage{amsfonts}
\usepackage{amssymb}
\usepackage{amsmath}
\usepackage{mathrsfs}
\input xypic
\usepackage{hyperref}

\newcommand{\be}{\begin{equation}}
\newcommand{\ee}{\end{equation}}
\newcommand{\bea}{\begin{eqnarray}}
\newcommand{\eea}{\end{eqnarray}}
\newcommand{\bean}{\begin{eqnarray*}}
\newcommand{\eean}{\end{eqnarray*}}
\newcommand{\brray}{\begin{array}}
\newcommand{\erray}{\end{array}}

\newtheorem{dfn}{Definition}[section]
\newtheorem{thm}[dfn]{Theorem}
\newtheorem{lmma}[dfn]{Lemma}
\newtheorem{ppsn}[dfn]{Proposition}
\newtheorem{crlre}[dfn]{Corollary}
\newtheorem{xmpl}[dfn]{Example}
\newtheorem{rmrk}[dfn]{Remark}
\newcommand{\bdfn}{\begin{dfn}\rm}
\newcommand{\bthm}{\begin{thm}}
\newcommand{\blmma}{\begin{lmma}}
\newcommand{\bppsn}{\begin{ppsn}}
\newcommand{\bcrlre}{\begin{crlre}}
\newcommand{\bxmpl}{\begin{xmpl}}
\newcommand{\brmrk}{\begin{rmrk}\rm}

\newcommand{\edfn}{\end{dfn}}
\newcommand{\ethm}{\end{thm}}
\newcommand{\elmma}{\end{lmma}}
\newcommand{\eppsn}{\end{ppsn}}
\newcommand{\ecrlre}{\end{crlre}}
\newcommand{\exmpl}{\end{xmpl}}
\newcommand{\ermrk}{\end{rmrk}}

\newcommand{\bbc}{\mathbb{C}}
\newcommand{\bbz}{\mathbb{Z}}

\newcommand{\bbn}{\mathbb{N}}

\newcommand{\bbt}{\mathbb{T}}

\newcommand{\cla}{\mathcal{A}}
\newcommand{\clb}{\mathcal{B}}

\newcommand{\clh}{\mathcal{H}}

\newcommand{\cls}{\mathcal{S}}

\title{The weak heat kernel asymptotic expansion and the quantum double
suspension}
\author{Partha Sarathi Chakraborty and S.Sundar}

\begin{document}
\maketitle
\begin{abstract}
In this paper we are concerned with the construction of a general principle that will allow us to produce regular spectral triples with finite and simple dimension spectrum.  We introduce the notion of weak heat kernel asymptotic expansion (WHKAE) property of a spectral triple and show that the weak heat kernel asymptotic expansion allows one to conclude that the spectral triple is regular with finite simple dimension spectrum. The usual heat kernel expansion implies this property. Finally we show that WHKAE is stable under quantum double suspension, a notion introduced in \cite{Hong-Sym1}. Therefore quantum double suspending compact Riemannian spin manifolds iteratively we get many examples of regular spectral triples with finite simple dimension spectrum. This covers all the odd dimensional quantum spheres. Our methods also apply to the case of noncommutative torus.
\end{abstract}
\section{Introduction}
Since it's inception index theorems play a central role in noncommutative geometry.
Here spaces are replaced by explicit K-cycles or finitely summable Fredholm modules. Through index pairing they pair naturally with K-theory. In the foundational paper Alain Connes introduced cyclic cohomology as a natural recipient of a Chern character homomorphism assigning cyclic cocycles to finitely summable Fredholm modules. Then the index pairing is computed by the pairing of cyclic cohomology with K-theory. But finitely summable Fredholm modules occur as those associated with spectral triples and it was desirable to have cyclic cocycles given directly in terms of spectral data, that will compute the index pairing. This was achieved by Connes and Moscovici in \cite{Con-Mos}. Let us briefly recall their local index formula henceforth to be abbreviated as LIF. One begins with a spectral triple, i.e., a Hilbert space ${\mathcal H}$, an involutive subalgebra $\mathcal A$ of the algebra of bounded operators on $\mathcal H$ and a self adjoint operator $D$ with compact resolvent. It is  further assumed that the commutators $[D,{\mathcal A}]$ give rise to bounded operators. Such a triple is finitely summable if ${|D|}^{-p}$ is trace class for some positive $p$. The spectral triple is said to be regular if both $\mathcal A$ and $[D,{\mathcal A}]$ are in the domains of  $\delta^n$ for all $n \ge 0$, where $\delta$ is the derivation $[|D|,\cdot]$. One says that the spectral triple has  dimension spectrum $\Sigma$, if for every element $b$ in the smallest algebra $\mathcal B$ containing $\mathcal A$, $[D,{\mathcal A}]$ and closed under the derivation $\delta$, the associated zeta function $\zeta_b(z)=Tr b {|D|}^{-z}$ a priori defined on the right half plane $\Re (z) > p$ admits a meromorphic extension to whole of complex plane with poles contained in $\Sigma$. If a spectral triple is regular and has discrete dimension spectrum then given any $n$-tuple of non-negative integers $k_1, k_2, \cdots, k_n$ one can consider multilinear functionals  $\phi_{{\bf k},n}$ defined by $ \phi_{{\bf k},n}(a_0,a_1,\cdots,a_n)=Res_{z=0} Tr a_0 {[D,a_1]}^{(k_1)} {[D,a_2]}^{(k_2)}\cdots {[D,a_n]}^{(k_n)}{|D|}^{-n-2|k|-z}$, where $T^{(r)}$ stands for the $r$-fold commutator $[D^2,[D^2,[\cdots[D^2 ,T]\cdots]$ and $|k|=k_1+\cdots+k_n$. In the local index formula  the components of the local Chern character in the $(b,B)$-bicomplex is expressed as a sum $\sum_{{\bf k}} c_{{\bf k},n} \phi_{{\bf k},n}$, where the summation is over all $n$-tuples of non negative integers and $c_{{\bf k},n}$'s are some universal constants independent of the particular spectral triple under consideration. Note that  remark II.1 in page 63 of \cite{Con-Mos} says that if we consider the Dirac operator associated with a closed Riemannian spin manifold then $\phi_{{\bf k},n}$'s are zero for $|k| \ne 0$. Therefore most of the terms in the local Chern character are visible in truely noncommutative cases and hence should be interpreted as a signature of noncommutativity. To have a better understanding of the contribution of these terms it is desirable to have examples where these terms  survive. Even though the Connes-Moscovici paper ended with the hope that in the foliations example these contributions will be visible it turned out that ``even in the case of codimension one foliations, the printed form of the cocycle takes around one hundred pages'' \cite{Con-Mos-Hopf}. So  the task of illustrating the LIF in simpler examples remained open. The first illustration was given by Connes in \cite{ConSU_q2}. This was extended to odd dimensional quantum spheres by Pal and Sundar in \cite{SunPal}. But to have a good grasp of the formula it is essential to have a systematic family of examples where one can verify the hypothesis of regularity and discreteness of the dimension spectrum. To our knowledge only Higson (\cite{Higson}) made an attempt to this effect. Here in this article we are also primarily concerned with that goal, namely developing a general procedure that will allow us to construct regular spectral triples with finite dimension spectrum from known examples. As in \cite{Higson} we also draw inspiration from the classical situation. We show that a hypothesis similar to, but weaker than the heat kernel expansion which we call weak heat kernel asymptotic expansion implies regularity and discreteness of the dimension spectrum. The usual heat kernel expansion implies the weak heat kernel expansion. More importantly we show that the weak heat kernel expansion is stable under quantum double suspension, a notion introduced in \cite{Hong-Sym1}. Therefore by iteratively quantum double suspending compact Riemannian manifolds we get examples of noncommutative geometries which are regular with finite dimension spectrum. We show noncommutative torus satisfies the weak heat kernel expansion and there by satisfies regularity and discreteness of it's dimension spectrum.

Organization of the paper is as follows. In section 2 we recall the basics of Mellin transform and asymptotic expansions. In the next section we introduce the weak heat kernel expansion property and show that this implies regularity and finiteness of the dimension spectrum. We also show that the usual heat kernel expansion implies the weak heat kernel expansion. In the next section we recall the notion of quantum double suspension and show that weak kernel expansion is stable under quantum double suspension. The weak heat kernel expansion of noncommutative torus is also established. In the final section we do a topological version of the theory relevant for applications in quantum homogeneous spaces. 
\section{Asymptotic expansions and the Mellin transform}
 First we recall a few basic facts about asymptotic expansions. Let
$\phi:(0,\infty)\to \mathbb{C}$ be a continuous function. We say that $\phi$ has
an asymptotic power series expansion near $0$ if there exists a sequence
$(a_{r})_{r=0}^{\infty}$ of complex numbers such that given $N$ there exists
$\epsilon,M >0$ such that if $t \in(0,\epsilon)$
\begin{displaymath}
 |\phi(t)-\sum_{r=0}^{N}a_{r}t^{r}| \leq M t^{N+1}.
\end{displaymath}
We write $\phi(t) \sim \sum_{0}^{\infty}a_{r}t^{r}$ as $t\to 0+$. Note that the
coefficients $(a_{r})$ are unique. For, \begin{equation}
                 a_{N} = \lim_{t\to
0+}\frac{\phi(t)-\sum_{r=0}^{N-1}a_{r}t^{r}}{t^{N}}.
               \end{equation}
If $\phi(t)\sim \sum_{r=0}^{\infty}a_{r}t^{r}$ as $t \to 0+$ then $\phi$ can be
extended continuously to $[0,\infty)$ simply by letting $\phi(0):=a_{0}$.

Let $X$ be a topological space and $F:[0,\infty)\times X \to \bbc$ be a
continuous function. Suppose that for every $x \in X$, the function $t\to
F(t,x)$ has an asymptotic expansion near $0$
\begin{equation}
\label{uniform}
 F(t,x) \sim \sum_{r=0}^{\infty}a_{r}(x)t^{r}.
\end{equation}
Let $x_{0} \in X$. We say that Expansion \ref{uniform} is uniform at $x_{0}$ if
given $N$ there exists an open set $U \subset [0,\infty)\times X$ containing
$(0,x_{0})$ and an $M>0$ such that for $(t,x)\in U$ one has  
\begin{displaymath}
 |F(t,x)-\sum_{r=0}^{N}a_{r}(x)t^{r}| \leq M t^{N+1}.
\end{displaymath}
 We say that Expansion \ref{uniform} is uniform if it is uniform at every point
of $X$.

\begin{ppsn}
 Let $X$ be a topological space and $F:[0,\infty)\times X \to \bbc$ be a
continuous function. Suppose that $F$ has a uniform asymptotic power series
expansion
\begin{displaymath}
 F(t,x) \sim \sum_{r=0}^{\infty}a_{r}(x)t^{r}.
\end{displaymath}
Then for every $r \geq 0$, the function $a_{r}$ is continuous.
\end{ppsn}
\textit{Proof.} It is enough to show that the function $a_{0}$ is continuous.
Let $x_{0} \in X$ be given. Since the expansion of $F$ is uniform at $x_{0}$, it
follows that there exists an open set $U$ containing $x_{0}$ and $\delta,M >0$
such that 
\begin{equation}
\label{inequality}
 |F(t,x)-a_{0}(x)| \leq M t \mbox{ for } t<\delta \mbox{ and } x\in U.
\end{equation}
Let $F_{n}(x):=F(\frac{1}{n},x)$. Then Equation \ref{inequality} says that
$F_{n}$ converges uniformly to $a_{0}$ on $U$. Hence $a_{0}$ is continuous on
$U$ and hence at $x_{0}$. This completes the proof. \hfill $\Box$ 

The following two lemmas are easy to prove and we leave the proof to the reader.
\begin{lmma}
\label{product}
 Let $X ,Y$ be topological spaces. Let $F:[0,\infty)\times X \to \bbc$ and
$G:[0,\infty)\times Y \to \bbc$ be continuous. Suppose that $F$ and $G$ has
uniform asymptotic power series expansion. Then the function $H:[0,\infty)\times
X \times Y \to \bbc$ defined by $H(t,x,y):=F(t,x)G(t,y)$ has uniform asymptotic
power series expansion.\\
Moreover if \begin{displaymath}
             F(t,x)\sim \sum_{r=0}^{\infty}a_{r}(x)t^{r} \mbox{ and } G(t,y)\sim
\sum_{r=0}^{\infty}b_{r}(y)t^{r},
            \end{displaymath}
then
\begin{displaymath}
  H(t,x,y) \sim \sum_{r=0}^{\infty}c_{r}(x,y)t^{r},
\end{displaymath}
where
\begin{displaymath}
 c_{r}(x,y):=\sum_{m+n=r}a_{m}(x)b_{n}(y).
\end{displaymath}
\end{lmma}

\begin{lmma}
\label{entire}
 Let $\phi:[1,\infty)\to \bbc$ be a continuous function. Suppose that for every
$N$, $$\sup_{t\in[1,\infty)}|t^{N}\phi(t)|< \infty.$$ Then the function $s
\mapsto \int_{1}^{\infty}\phi(t)t^{s-1}dt$ is entire.
\end{lmma}
\subsection{The Mellin transform}
    In this section we recall the definition of the Mellin transform of a
function defined on $(0,\infty)$ and analyse the relationship between the
asymptotic expansion of a function and the meromorphic continuation of its
Mellin transform. Let us introduce some notations. We say that a function
$\phi:(0,\infty) \to \mathbb{C}$ is of rapid decay near infinity if for every
$N>0$, $\sup_{t\in[1,\infty)}|t^{N}\phi(t)|$ is finite. We let
$\mathcal{M}_{\infty}$ to be the set of continuous complex valued functions on
$(0,\infty)$ which has rapid decay near infinity.
 For $p \in \mathbb{R}$, we let  
\begin{align*}
 \mathcal{M}_{p}((0,1]):&=\{\phi:(0,1]\to \mathbb{C}: \text{ $\phi$ is
continuous and } ~\sup_{t \in (0,1]}t^{p}|\phi(t)|<\infty\},\\
\mathcal{M}_{p}:&=\{\phi \in \mathcal{M}_{\infty}: \phi|_{(0,1]} \in
\mathcal{M}_{p}((0,1])\} .
\end{align*}
Note that if $p \leq q$ then $\mathcal{M}_{p} \subset \mathcal{M}_{q}$ and
$M_{p}((0,1]) \subset M_{q}((0,1])$.
\begin{dfn}
 Let $\phi:(0,\infty) \to \mathbb{C}$ be a continuous function. Suppose that
$\phi \in \mathcal{M}_{p}$ for some $p$.
 Then the Mellin transform of $\phi$, denoted $M\phi$, is defined as follows:
For $Re(s)>p$,
 \begin{equation*}
  M\phi(s):=\int_{0}^{\infty}\phi(t)t^{s-1}dt.
 \end{equation*}
\end{dfn}
 One can show that if $\phi \in \mathcal{M}_{p}$ then $M\phi$ is analytic on the
right half plane $Re(s)>p+2$. Also if $\phi \in \mathcal{M}_{p}((0,1])$ then $s
\mapsto \int_{0}^{1}\phi(t)t^{s-1}$ is analytic on $Re(s)>p+2$. 

For $a<b$ and
$K>0$, let $H_{a,b,K}:=\{\sigma+it: ~a\leq \sigma \leq
b,~|t|>K \}$.
\begin{dfn}
 Let $F$ be a meromorphic function on the entire complex plane with simple poles
lying inside the set of integers. We say that $F$ has decay of order $r \in
\bbn$ along the vertical strips if  the function 
$s\mapsto s^{r}F(s)$ is bounded on $H_{a,b,K}$ for every $a<b$ and $K>0$. We say
that $F$ is of rapid decay along the vertical strips if $F$ has decay of order
$r$ for every $r \in \bbn$.

\end{dfn}

\begin{ppsn}
\label{MTT}
 Let $\phi:(0,\infty)\to \mathbb{C}$ be a continuous function of rapid decay.
Assume that $\phi(t) \sim \sum_{0}^{\infty}a_{r}t^{r}$ as $t\to 0+$.  Then we
have the following. 
\begin{itemize}
 \item[(1)] The function $\phi \in \mathcal{M}_{0}$.
 \item[(2)] The Mellin transform $M\phi$ of $\phi$ extends to a meromorphic
function to the whole of complex plane with simple poles in the set of negative
integers. $\{0,-1,-2,-3,\cdots\}$
 \item[(3)] The residue of $M\phi$ at $s=-r$ is given by $Res_{s=-r}M\phi(s)=
a_{r}$.
\item[(4)] The meromorphic continuation of the Mellin transform $M\phi$ has
decay of order $0$ along the vertical strips.
\end{itemize}
\end{ppsn}
\textit{Proof.} By definition it follows that $\phi \in \mathcal{M}_0$. Since
$\phi$ has rapid decay at infinity, by lemma \ref{entire}, it follows that the
function $s \mapsto \int_{1}^{\infty}\phi(t)t^{s-1}dt$ is entire. Thus modulo a
holomorphic function $M\phi(s) \equiv \int_{0}^{1}\phi(t)t^{s-1}$. For $N \in
\bbn$, Let $R_{N}(t):=\phi(t)-\sum_{r=0}^{N}a_{r}t^{r}$. Thus modulo a
holomorphic function, we have
\[
 M\phi(s)\equiv \sum_{r=0}\frac{a_{r}}{s+r}+\int_{0}^{1}R_{N}(t)t^{s-1}dt
\]
As $R_{N} \in \mathcal{M}_{-(N+1)}((0,1])$ the function $s\mapsto
\int_{0}^{1}R_{N}(t)t^{s-1} dt$ is holomorphic on $Re(s)>-N+1$. Thus on
$Re(s)>-N+1$, modulo a holomorphic function, one has
\begin{equation}
\label{asymp}
 M\phi(s)\equiv \sum_{r=0}^{N}\frac{a_{r}}{s+r}.
\end{equation}
This shows that $M\phi$ admits a meromorphic continuation to the whole of
complex plane and has simple poles lying in the set of negative integers
$\{0,-1,-2,\cdots,\}$. Also $(3)$ follows from Equation \ref{asymp}.

Let $a<b$ and $K>0$ be given. Choose $N \in \mathbb{N}$ such that $N+a>0$. Then
one has
\[
 M\phi(s)=\sum_{r=0}^{N}\frac{a_{r}}{s+r}+\int_{0}^{1}R_{N}(t)t^{s-1}dt+
\int_{1}^{\infty}\phi(t)t^{s-1}dt.
\]
As the function $s \mapsto \frac{1}{s+r}$ is bounded for every $r\geq 0$ on
$H_{a,b,K}$, it is enough to show that the functions
$\psi(s):=\int_{0}^{1}R_{N}(t)t^{s-1}dt$ and
$\chi(s):=\int_{1}^{\infty}\phi(t)t^{s-1}dt$ are bounded on $H_{a,b,K}$.

By definition of the asymptotic expansion, it follows that there exists an $M>0$
such that $|R_{N}(t)|\leq Mt^{N+1}$. Hence for $s:=\sigma+it \in H_{a,b,K}$, 
\[
|\psi(s)| \leq  \frac{M}{\sigma+N+1} \leq  \frac{M}{a+N+1} \leq M .
\]
 Thus $\psi$ is bounded on $H_{a,b,K}$.

Now for $s:=\sigma+it \in H_{a,b,K}$, we have 
\[
|\chi(s)| \leq \int_{1}^{\infty}|\phi(t)|t^{\sigma-1}dt \leq
\int_{1}^{\infty}|\phi(t)|t^{b-1}dt .
\]
 Since $\phi$ is of rapid decay, the integral 
$\int_{1}^{\infty}|\phi(t)|t^{a-1}dt$ is finite. Hence $\chi$ is bounded on
$H_{a,b,K}$. This completes the proof. \hfill $\Box$      \begin{crlre}
              Let $\phi:(0,\infty)\mapsto \bbc$ be a smooth function. Assume
that for every $n$, the n$^{th}$ derivative $\phi^{(n)}$ has rapid decay at
infinity and admits an asymptotic power series expansion near $0$.              
                     
\begin{enumerate}
\item[(1)] For every $n$, the Mellin transform $M\phi^{(n)}$ of $\phi^{(n)}$
extends to a meromorphic function to the whole of complex plane with simple
poles in the set of negative integers $\{0,-1,-2,-3,\cdots\}$.
\item[(2)] The meromorphic continuation of the Mellin transform $M\phi$ is of
rapid decay along the vertical strips.
\end{enumerate}
\end{crlre}
\textit{Proof.} (1) follows from Proposition \ref{MTT}. To prove (2), observe
that $M\phi^{'}(s+1)= -sM\phi(s)$. For $Re(s)\gg 0$,
\begin{align*}
 M\phi^{'}(s+1):&= \int_{0}^{\infty}\phi^{'}(t)t^{s}dt\\
                &= -\int_{0}^{\infty}s\phi(t)t^{s-1}dt \textrm{ (~follows from
integration by parts)}\\
                &= -sM\phi(s).
\end{align*}
As $M\phi^{'}$ and $M\phi$ are meromorphic, it follows that $M\phi^{'}(s+1)=
-sM\phi(s)$. Now a repeated application of this equation gives
\begin{equation}
 \label{recursion}
M\phi(s):=(-1)^{n}\frac{M\phi^{(n)}(s+n)}{s(s+1)\cdots (s+n-1)}.
\end{equation}
 Now let $a<b, K>0$ and $r \in \bbn$ be given. Now $(3)$ of Proposition
\ref{MTT} applied to $\phi^{(r)}$ , together with Equation \ref{recursion},
implies that the function $s\mapsto s^{r}M\phi(s)$ is bounded on $H_{a,b,K}$.
This completes the proof. \hfill $\Box$

The following proposition shows how to pass from the decay properties of the
Mellin transform of a function to the asymptotic expansion property of the
function.
\begin{ppsn}
\label{IMTT}
 Let $\phi \in \mathcal{M}_{p}$ for some $p$. Assume that the Mellin transform
$M\phi$ is meromorphic on the entire complex plane with  poles lying in the set
of negative integers $\{0,-1,-2,\cdots\}$. Suppose that the meromorphic
continuation of the Mellin transform $M\phi$ is of rapid decay along the
vertical strips. Then the function $\phi$ has an asymptotic expansion near $0$. 

Moreover if $a_{r}:= Res_{s=-r}M\phi(s)$ then $\phi(t)\sim
\sum_{r=0}^{\infty}a_{r}t^{r}$ near $0$.
\end{ppsn}
\textit{Proof.} The proof is a simple application of the inverse Mellin
transform. Let $M\gg0$. Then one has the following inversion formula.
\begin{equation*}
 \phi(t)=\int_{M-i\infty}^{M+\infty}M\phi(s)t^{-s}ds
\end{equation*}
Define $F_{t}(s):=M\phi(s)t^{-s}$. Suppose $N \in \bbn$ be given. Let $\sigma
\in (-N-1,-N)$ be given. For every $A>0$, by Cauchy's integral formula, we have 
\begin{equation}
\label{cauchy}
 \int_{M-iA}^{M+iA}F_{t}(s)ds+\int_{M+iA}^{\sigma+iA}F_{t}(s)ds+\int_{\sigma+iA}
^{\sigma-iA}F_{t}(s)ds+\int_{\sigma-iA}^{M-iA}F_{t}(s)ds=
\sum_{r=0}^{N}Res_{s=-r}F_{t}(s)
\end{equation}
For a fixed $t$, $F_{t}$ has rapid decay along the vertical strips. Thus when $A
\to \infty$ the second and fourth integrals in Equation \ref{cauchy} vanishes
and we obtain the following equation
\begin{equation}
\label{estimate0}
 \phi(t)-\sum_{r=0}^{N}a_{r}t^{r} = \int_{\sigma-i\infty}^{\sigma+i\infty}
M\phi(s)t^{-s}ds
\end{equation}
But $M\phi(\sigma+it)$ has rapid decay in $t$. Let
$M_{\sigma}:=\int_{-\infty}^{\infty}|M\phi(\sigma+it)|$. Then Equation
\ref{estimate0} implies that 
\[
|\phi(t)-\sum_{r=0}^{N}a_{r}t^{r}| \leq M_{\sigma}t^{-\sigma} \leq M_{\sigma}
t^{N} \text{~for }t \leq 1
\]
 Thus we have shown that for every $N$,
$R_{N}(t):=\phi(t)-\sum_{r=0}^{N}a_{r}t^{r}= O(t^{N})$ as $t \to 0$ and hence
$R_{N-1}(t)=R_{N}(t)+a_{N}t^{N}=O(t^{N})$  as $t \to 0$.
 This completes the proof. \hfill $\Box$

\section{The weak heat kernel asymptotic expansion property and the dimension
spectrum of  spectral triples}
        In this section we consider a property of spectral triples which we call
the weak heat kernel asymptotic expansion property. We show that a spectral
triple having the weak heat kernel asymptotic expansion property is regular and
has finite dimension spectrum lying in the set of postive integers.
\begin{dfn}
\label{heat kernel}
Let $(\mathcal{A},\mathcal{H},D)$ be a $p+$ summable spectral
triple for a $C^{*}$ algebra A where $\mathcal{A}$ is a dense $*$
subalgebra of $A$. We say that the spectral triple
$(\mathcal{A}, \mathcal{H}, D)$ has the weak heat kernel asymptotic
expansion property if there exists a $*$ subalgebra $\mathcal{B}
\subset B(\mathcal{H})$ such that 
\begin{itemize}
\item[(1)] The algebra $\mathcal{B}$ contains $\mathcal{A}$,
\item[(2)] The unbounded derivations $\delta:=[|D|,.]$  leaves
$\mathcal{B}$ invariant. Also the unbounded derivation $d:=[D,.]$ maps
$\mathcal{A}$ into $\mathcal{B}$.
\item[(3)] The algebra $\mathcal{B}$ is invariant under the left multiplication
by $F$ where $F:=sign(D)$.
\item[(4)] For every $b \in \mathcal{B}$, the function $\tau_{p,b}:(0,\infty)
\mapsto \mathbb{C}$ defined by
$\tau_{p,b}(t)=t^{p}Tr(be^{-t|D|})$ has an  asymptotic power series
expansion.
\end{itemize}
\end{dfn}
If the algebra $\mathcal{A}$ is unital and the representation of $\cla$ on
$\clh$ is unital then $(3)$ can be replaced by the condition $F \in
\mathcal{B}$. The next
proposition proves that an odd spectral triple
that has the heat kernel asymptotic expansion property is regular and has simple
dimension spectrum.
\begin{ppsn}
\label{dimension spectrum}
Let $(\mathcal{A},\mathcal{H}, D)$ be a $p+$ summable 
spectral triple which has the weak heat kernel asymptotic expansion property.
Then
the spectral triple $(\mathcal{A}, \mathcal{H}, D)$ is regular and
has finite simple dimension spectrum. Moreover the dimension spectrum is
contained in $\{1,2,\cdots,p\}$.
\end{ppsn}
\textit{Proof.} Let $\mathcal{B} \subset B(\mathcal{H})$ be a $*$
algebra for which $(1)-(4)$ of Definition \ref{heat kernel} is satisfied. The
fact that $\mathcal{B}$ satisfies $(1)$ and $(2)$ imples that  the
spectral triple
$(\mathcal{A},\mathcal{H},D)$ is regular. First we assume that $D$ is
invertible. Let $b \in \mathcal{B}$ be given. 

Since $|D|^{-q}$ is trace class for  $q>p$, it follows that for every $N
>p$ there exists an $M>0$ such that $Tr(e^{-t|D|}) \leq M
t^{-N}Tr(|D|^{-N})$. Now for $1 \leq t<\infty$ and $N \geq p$ one has
\begin{align*}
 |Tr(be^{-t|D|})| ~& \leq ~\|b\|Tr(e^{-t|D|}) \\
                  ~& \leq ~\|b\|M t^{-N}Tr(|D|^{-N})
\end{align*}
Thus the function $t\mapsto Tr(be^{-t|D|})$ is of  rapid decay near infinity.
Now observe that for $Re(s)\gg 0$
\begin{equation}
\label{MT of heat kernel}
Tr(b|D|^{-s})= \frac{1}{\Gamma(s)}\int_{0}^{\infty}Tr(be^{-t|D|})t^{s-1}dt\\
\end{equation} 

 By assumption the function $\phi(t):= t^{p}Tr(be^{-t|D|})$ has an asymptotic
power series expansion near $0$. By Equation \ref{MT of heat kernel}, it follows
that $M\phi(s)=\Gamma(s+p)Tr(b|D|^{-s-p})$. Now  Propostition \ref{MTT} implies
that the function  $s \mapsto \Gamma(s)Tr(b|D|^{-s})$ is meromorphic with simple
poles lying inside $\{n \in \mathbb{Z}: n \leq p\}$. As $\frac{1}{\Gamma(s)}$ is
entire and has simple zeros at $\{k: k \leq 0\}$, it follows that the function
$s \to Tr(b|D|^{-s})$ is meromorphic and has simple poles with poles lying in
$\{1,2,\cdots,p\}$.

Suppose $D$ is not invertible. Let $P$ denote the projection onto the kernel of
$D$ which is finite dimensional. Let $D^{'}:=D+P$ and $b$ be an element in
$\mathcal{B}^{\infty}$. Now note that 
\begin{displaymath}
Tr(be^{-t|D^{'}|})=Tr(PbP)e^{-t} + Tr(be^{-t|D|}) 
\end{displaymath}
 Hence the function $t \to t^{p}Tr(be^{-t|D^{'}|})$ has an asymptotic power
series expansion. Thus the function $s \to Tr(b|D^{'}|^{-s})$ is meromorphic
with simple poles lying in $\{1,2,\cdots,p\}$. Observe that for $Re(s)\gg 0$,
$Tr(b|D^{'}|^{-s})=Tr(b|D|^{-s})$.  Hence the function $s \to Tr(b|D|^{-s})$ is
meromorphic with simple poles lying in $\{1,2,\cdots,p\}$. This completes the
proof.
\hfill $\Box$ 
\begin{rmrk}
\label{residue}
 If $Tr(be^{-t|D|})\sim \sum_{r=-p}^{\infty}a_{r}(b)t^{r}$ then (3) of
Proposition \ref{MTT} imples that 
\begin{align*}
Res_{z=k}Tr(b|D|^{-z})&= \frac{1}{k!}a_{-k}(b) \textrm{~for~} 1 \leq k \leq p\\
Tr(b|D|^{-z})_{z=0}&=a_{0}(b) 
\end{align*}
\end{rmrk}
\begin{rmrk}
 Let $(\mathcal{A}, \mathcal{H}, D)$ be a spectral triple which has the
weak heat kernel asymptotic expansion property. Then the dimension spectrum
$\Sigma$
is finite and lies in the set of positive integers. We call the greatest element
in the dimension spectrum as the dimension of the spectral triple
$(\mathcal{A},\mathcal{H},D)$. If $\Sigma$ is empty we set the
dimension to be $0$.
\end{rmrk}
Now in the next proposition we show that the usual heat kernel asymptotic
expanansion implies the weak heat kernel asymptotic expansion. 

\begin{ppsn}
\label{implication}
 Let $(\mathcal{A}, \mathcal{H}, D)$ be a $p+$ summable spectral
triple for a $C^{*}$ algebra A. 
Suppose that $\mathcal{B}$ is a $*$ subalgebra of $B(\clh)$
satisfying  $(1)-(4)$ of Definition \ref{heat kernel}. Assume that for every
$b \in \mathcal{B}$, the function $\sigma_{p,b}:(0,\infty)
\to \mathbb{C}$ defined by $\sigma_{p,b}(t):=t^{p}Tr(be^{-t^{2}D^{2}})$ has an
 asymptotic power series expansion. 

Then for every $b \in \mathcal{B}$, the function $\tau_{p,b}:t \mapsto
t^{p}Tr(be^{-t|D|})$ has an asymptotic power series expansion.
\end{ppsn}
\textit{Proof.} It is enough to consider the case where $D$ is invertible. Let
$b \in \mathcal{B}$ be given. Let $\psi$ denotes the Mellin transform
of the function $t \mapsto Tr(be^{-t^{2}D^2})$ and $\chi$ denote the Mellin
transform of the function $t\mapsto Tr(be^{-t|D|})$. Then a simple change of
variables shows that $\psi(s)=\frac{\Gamma(\frac{s}{2})}{2}Tr(b|D|^{-s})$. But
then $\chi(s)=\Gamma(s)Tr(b|D|^{-s})$. Thus we obtain the equation 
\[
\chi(s)=\frac{2\Gamma(s)}{\Gamma(\frac{s}{2})}\psi(s)
\]
 But we have following duplication
formula for the gamma function 
\[
 \Gamma(s)\Gamma(s+\frac{1}{2})= 2^{1-2s}\sqrt{\pi}\Gamma(2s)
\]
Hence one has 
\[
\chi(s)=\frac{1}{\sqrt{\pi}}2^{s}\Gamma(\frac{s+1}{2})\psi(s) 
\]
Now Proposition \ref{MTT} implies that $\psi$ has  decay of order $0$ along the
vertical strips and has simple poles lying inside $\{n \in \mathbb{Z}:n \leq
p\}$. Since the gamma function has rapid decay along the vertical strips, it
follows that $\chi$ has rapid decay along the vertical strips and has poles
lying in $\{n \in \mathbb{Z}:n \leq p\}$. If $\tilde{\chi}$ denotes the Mellin
transform of $\tau_{p}(.,b)$ then $\tilde{\chi}(s)=\chi(s+p)$. Hence
$\tilde{\chi}$ has rapid decay along the vertical strips and has poles lying in
the set of negative integers. Now Proposition \ref{IMTT} implies that the map $t
\to t^{p}Tr(be^{-t|D|})$ has an asymptotic power series expansion near $0$. This
completes the proof. \hfill $\Box$
\section{Stability of the weak heat kernel expansion property and the quantum
double suspension}

Let us recall the  definition of the quantum double suspension of a
unital $C^{*}$ algebra . The quantum double suspension is first defined in
\cite{Hong-Sym1} and our equivalent definition is as in \cite{Hong-Sym2}. Let us
fix some notations. We denote the left shift on $\ell^{2}(\mathbb{N})$ by $S$
which is defined on the standard orthonormal basis $(e_{n})$ as $Se_{n}=e_{n-1}$
and $p$ denote the projection $|e_{0}\rangle \langle e_{0}|$. The number
operator on $\ell^{2}(\bbn)$ is denoted by $N$ and defined as $Ne_{n}:=ne_{n}$.
We denote the
$C^{*}$ algebra generated by $S$ in $B(\ell^{2}(\mathbb{N}))$ by $\mathscr{T}$
which is the toeplitz algebra. Note that $SS^{*}=1$ and $p=1-S^{*}S$. Let
$\sigma:\mathscr{T} \to C(\mathbb{T})$ be the symbol map which sends $S$ to the
generating unitary $z$. Then one has the following exact sequence
\begin{displaymath}
 0 \to \mathcal{K} \to \mathscr{T} \stackrel{\sigma} \to C(\mathbb{T}) \to 0
\end{displaymath}

\begin{dfn}
 Let $A$ be a unital $C^{*}$ algebra. Then the quantum double suspension of $A$
denoted $\Sigma^{2}(A)$ is the $C^{*}$ algebra generated by $A\otimes p$ and
$1\otimes S$ in $A\otimes \mathscr{T}$.
\end{dfn}
Let $A$ be a unital $C^{*}$ algebra. One has the following exact sequence.
\begin{displaymath}
 0 \to A\otimes \mathcal{K}(\ell^{2}(\mathbb{N})) \to \Sigma^{2}(A)
\stackrel{\rho} \to C(\mathbb{T}) \to 0
\end{displaymath}
where $\rho$ is just the restriction of $1\otimes \sigma$ to $\Sigma^{2}(A)$. 
\begin{rmrk}
 It can be easily shown that $\Sigma^{2}(C(\mathbb{T}))=C(SU_{q}(2))$ and more
generally one can show that $\Sigma^{2}(C(S_{q}^{2n-1}))=C(S_{q}^{2n+1})$. We
refer to \cite{Hong-Sym1} or Lemma 3.2 of \cite{SunPal} for the proof. Thus the
odd dimensional quantum spheres can be obtained from the circle $\mathbb{T}$ by
applying the quantum double suspension recursively.
\end{rmrk}

Let
$\mathcal{A}$ be a dense $*$ subalgebra of a $C^{*}$ algebra $A$. 
Define \[  
\Sigma^{2}_{alg}(\mathcal{A}):= span\{ a \otimes k, 1\otimes S^{n}, 1 \otimes
S^{*m} : a \in \mathcal{A}, k \in \mathcal{S}(\ell^{2}(\bbn)), n,m \geq 0\}
                  \]
where $\mathcal{S}(\ell^{2}(\bbn)):=\{(a_{mn}):
\sum_{m,n}(1+m+n)^{p}|a_{mn}|< \infty \text{~for every ~} p \}$.

Then $\Sigma^{2}_{alg}(\mathcal{A})$ is just the $*$ algebra generated by
$\mathcal{A}\otimes_{alg} \mathcal{S}(\ell^{2}(\bbn)$ and $1\otimes S$. Clearly
$\Sigma^{2}_{alg}(\mathcal{A})$ is a dense subalgebra of $\Sigma^{2}(A)$.

\begin{dfn}
 Let $(\mathcal{A},\clh, D)$ be a spectral triple and denote the sign of the
operator $D$ by $F$. Then the spectral triple
$(\Sigma^{2}_{alg}(\mathcal{A}),\clh \otimes \ell^{2}(\bbn),
\Sigma^{2}(D):=((F\otimes 1)(|D|\otimes 1+ 1\otimes N))$  is called the quantum
double suspension of the spectral triple $(\mathcal{A},\clh,D)$.
\end{dfn}
\subsection{ Stability of the weak heat kernel expansion}
We consider the stability of the weak heat kernel expansion under quantum
double suspension. First observe that the following are easily verifiable.
\begin{enumerate}
\item[(1)] The spectral triple $(\mathcal{S}(\ell^{2}(\bbn)),\ell^{2}(\bbn),N)$
has the weak heat kernel asymptotic expansion with dimension $0$.
 \item[(2)] Let $(\mathcal{A}_{i}, \mathcal{H}_{i},D_{i})$ be a spectral
triple with the weak heat kernel asymptotic expansion property with dimension
$p_{i}$
for $1 \leq i \leq n$. Then the spectral triple
$(\oplus_{i=1}^{n}\mathcal{A}_{i},\oplus_{i=1}^{n}\mathcal{H}_{i},
\oplus_{i=1}^{n}D_{i})$ has the weak heat kernel expansion property with
dimension
$p:=max\{p_{i}:1 \leq i \leq n\}$.
\item[(3)] If $(\mathcal{A},\mathcal{H},D)$ is a spectral triple with
the weak heat kernel asymptotic expansion property and has dimension $p$ then
$(\mathcal{A},\mathcal{H},|D|)$ also has the weak heat kernel asymptotic
expansion with the same dimension $p$.
\item[(4)] Let $(\mathcal{A},\mathcal{H},D)$ be a spectral triple with
the weak heat kernel asymptotic expansion property with dimension $p$.
Then the amplification $(\mathcal{A}\otimes
1,\mathcal{H}\otimes \ell^{2}(\mathbb{N}),|D|\otimes 1 + 1\otimes N)$ also has
the
asymptotic expansion property with dimension  $p+1$. 
\end{enumerate}

We start by proving the stability of the weak heat kernel expansion under
tensoring by compacts.
\begin{ppsn}
\label{tensoring by smooth compacts}
 Let $(\mathcal{A},\mathcal{H},D)$ be a spectral triple with the weak heat
kernel asymptotic expansion property of dimension $p$. Then
$(\mathcal{A}\otimes_{alg}
\mathcal{S}(\ell^{2}(\mathbb{N})),\mathcal{H}\otimes
\ell^{2}(\mathbb{N}),D_{0}:=(F\otimes 1)(|D|\otimes 1+1\otimes N))$ also has
the weak
heat kernel asymptotic expansion property with dimension $p$.
\end{ppsn}
\textit{Proof.} Let $\mathcal{B} \subset B(\mathcal{H})$ be a $*$
subalgebra for which $(1)-(4)$ of Definition \ref{heat kernel} are satisfied. We
denote $\mathcal{B}\otimes_{alg}\mathcal{S}(\ell^{2}(\mathbb{N}))$ by
$\mathcal{B}_{0}$. We show that $\mathcal{B}_{0}$ satisfies
$(1)-(4)$ of Definition \ref{heat kernel}. Clearly $(1)$ holds.

 We denote the unbounded derivation $[|D_{0}|,.]$,$[|D|,.]$ and $[N,.]$ by
$\delta_{D_{0}},\delta_{D}$ and $\delta_{N}$ respectively. By assumption
$\delta_{D}$ leaves $\mathcal{B}$ invariant.  Clearly
$\mathcal{B}\otimes_{alg}
\mathcal{S}(\ell^{2}(\mathbb{N}))$ is contained in the domain of
$\delta_{D_{0}}$ and
$\delta_{D_{0}}=\delta_{D}\otimes 1+1\otimes \delta_{N}$ on
$\mathcal{B}\otimes_{alg}
\mathcal{S}(\ell^{2}(\mathbb{N}))$.  Similarly one can show that the unbounded
derivation $[D_{0},.]$
maps $\cla\otimes_{alg}S(\ell^{2}(\bbn))$ into $\clb_{0}$ invariant.  

As $F_{0}:=sign(D_{0})=F\otimes 1$, $(3)$ is clear. Now $(4)$ follows
from Lemma\ref{product} and the equality $t^{p}Tr((b\otimes
k)e^{-t|D_{0}|}= t^{p}Tr(be^{-t|D|})Tr(ke^{-tN})$. This completes the proof.
\hfill $\Box$
    
Now we consider the stability of the heat kernel asymptotic expansion under the
double suspension. 

\begin{ppsn}
\label{suspension}
 Let $(\mathcal{A},\mathcal{H},D)$ be a spectral triple with the weak heat
kernel asymptotic expansion property of dimension $p$. Assume that the algebra
$\mathcal{A}$ is unital and the representation on $\mathcal{H}$ is
unital. Then the spectral triple $
(\Sigma^{2}(\mathcal{A}),\mathcal{H}\otimes
\ell^{2}(\mathbb{N}),\Sigma^{2}(D))$ also has the weak
heat kernel asymptotic expansion property with dimension $p+1$.
\end{ppsn}
\textit{Proof.} We denote $\Sigma^{2}(D)$ by $D_{0}$. Let
$\mathcal{B}$ be a $*$ subalgebra of $B(\mathcal{H})$
for which $(1)-(4)$ of Definition \ref{heat kernel} are satisfied. For
$f=\sum_{n}\lambda_{n}z^{n} \in C^{\infty}(\mathbb{T})$ , we let
$\sigma(f):= \sum_{n\geq 0}\lambda_{n}S^{n}+\sum_{n>0}\lambda_{-n}S^{*n}$. 
We denote the projection $\frac{1+F}{2}$ by $P$.  We let
$\mathcal{B}_{0}$ to denote
the algebra $\mathcal{B}\otimes_{alg}\cls(\ell^{2}(\bbn))$ as in
Proposition \ref{tensoring by smooth compacts}. As in Proposition \ref{tensoring
by smooth compacts}, we let $\delta_{D_{0}},\delta_{D},\delta_{N}$ to denote the
unbounded derivations $[|D_{0}|,.],[|D|,.]$ and $[N,.]$ respectively. Define
\begin{displaymath}
 \tilde{\mathcal{B}}:=\{b+P\otimes \sigma(f)+(1-P)\otimes \sigma(g): ~b
\in \mathcal{B}_{0},f,g \in C^{\infty}(\mathbb{T})\}
\end{displaymath}
Now it is clear that $\tilde{\mathcal{B}}$ satisfies $(1)$ of Definition
\ref{heat kernel}.

We have already shown in Proposition \ref{tensoring by smooth compacts} that
$\mathcal{B}_{0}$ is closed under $\delta_{D_{0}}$ and $d_{0}:=[D_{0},.]$ maps
$\mathcal{A}\otimes S(\ell^{2}(\bbn))$ into $\mathcal{B}_{0}$.
Now note that 
\begin{align*}
\delta_{D_{0}}(P\otimes \sigma(f))&=P\otimes \sigma(if^{'})\\
\delta_{D_{0}}((1-P)\otimes \sigma(g))&=(1-P)\otimes \sigma(ig^{'}) \\
[D_{0},P\otimes \sigma(f)]&= P \otimes \sigma(if^{'}) \\
[D_{0},(1-P)\otimes \sigma(g)]&=-(1-P)\otimes \sigma(ig^{'})
\end{align*}
Thus it follows that $\delta_{D_{0}}$ leaves $\tilde{\mathcal{B}}$ invariant
and $d_{0}:=[D_{0},.]$ maps $\Sigma_{2}(\cla)$ into
$\tilde{\clb}$. 

Since $F_{0}:=sign(D_{0})=F\otimes 1$, it follows from definition that $F_{0}
\in \tilde{\clb}$. Now we show that
$\tilde{\clb}$ satisfies $(4)$.

We have already shown in Proposition \ref{tensoring by smooth compacts} that
 given $b \in \clb_{0}$, the function $\tau_{p,b}(t)=t^{p}Tr(be^{-t|D_{0}|})$
has an asymptotic expansion. Hence the function $\tau_{p+1,b}$ has an
asymptotic expansion for every $b \in \clb_{0}$. Now note that 
\begin{eqnarray}
\label{asy}
 \tau_{p+1,P\otimes \sigma(f)}(t)&=& (\int f(\theta)d\theta)
t^{p}Tr(Pe^{-t|D|})tTr(e^{-tN})\\
\label{asy1}
 \tau_{p+1, (1-P)\otimes \sigma(g)}(t)&=&(\int
g(\theta)d\theta)t^{p}Tr((1-P)e^{-t|D|})tTr(e^{-tN})) 
\end{eqnarray}
Now recall that we have assumed that $\mathcal{A}$ is unital and hence
$P \in \mathcal{B}$. Hence $t^{p}Tr(xe^{-t|D|})$ has an asymptotic
power series expansion for $x\in \{P,1-P\}$. Thus $tTr(e^{-tN})$ has an
asymptotic power series expansion. From Equation \ref{asy}, Equation
\ref{asy1} and from  the earlier observation that $\tau_{p+1,b}$
 has an asymptotic power series expansion for $b \in \mathcal{B}_{0}$,
it follows that for every $b \in \tilde{\clb}$, the function $\tau_{p+1,b}$ has
an asymptotic power series expansion. This completes the proof.    \hfill $\Box$
\subsection{Higson's differential pair and the heat kernel expansion}
 
Now we discuss some examples of spectral triples which satisfy the weak heat
kernel asymptotic expansion property. In particular we discuss the spectral
triple associated to noncommutative torus and the classical spectral triple
associated to a spin manifold. Let us recall Higson's notion of a differential
pair as defined in \cite{Higson_localindex}. 

Consider  a Hilbert space $\clh$ and a positive, selfadjoint and
 an unbounded $\Delta$ on $\clh$. We assume that $\Delta$ has compact
resolvent.  For $k \in \bbn$, we let $\clh_{k}$ be the domain of the
operator $\Delta^{\frac{k}{2}}$. The vector space $\clh_{k}$ is given a Hilbert
space structure by identifying $\clh_{k}$ with the graph of the operator
$\Delta^{\frac{k}{2}}$. Denote the intersection $\bigcap_{k}\clh_{k}$ by
$\clh_{\infty}$.  An operator $T:\clh_{\infty}\to
\clh_{\infty}$ is said to be of analytic order  $\leq m$ where $m \in \bbz$
if $T$ extends to a bounded operator from $\clh_{k+m} \to \clh_{k}$ for every
$k$. We say an operator $T$ on $\clh_{\infty}$ has analytic order $-\infty$ if
$T$ has analytic order less than $-m$ for every $m>0$.  The following definition
is due to Higson. ( Refer \cite{Higson_localindex})

\begin{dfn}
 Let $\Delta$ be a positive, unbounded, selfadjoint operator on a Hilbert space
$\clh$ with compact resolvent. Suppose that $\mathcal{D}:=\bigcup_{p \geq
0}\mathcal{D}_{q}$ is a filtered algebra of operators on $\clh_{\infty}$. The
pair $(\mathcal{D},\Delta)$ is called a differential pair if the following
conditions hold.
\begin{enumerate}
 \item The algebra $\mathcal{D}$ is invariant under the derivation $T \to
[\Delta,T]$.
 \item If $X \in \mathcal{D}_{q}$, then $[\Delta,X] \in \mathcal{D}_{q+1}$.
 \item If $X \in \mathcal{D}_{q}$, then the analytic order of $X \leq q$.
\end{enumerate}
\end{dfn}

Now let us recall Higson's definition of pseudodifferential operators.
\begin{dfn}

 Let $(\mathcal{D},\Delta)$ be a differential pair. We denote the orthogonal
projection onto the kernel of $\Delta$ by $P$. Then $P$ is of finite rank as
$\Delta$ has compact resolvent. Let $\Delta_{1}:=\Delta+P$. Then $\Delta_{1}$
is invertible. 

 A  linear operator $T$ on
$\clh_{\infty}$ is called a basic pseudodifferential operator of order $ \leq k$
is  for every $\ell$ there exists $m$ and $X \in \mathcal{D}_{m+k}$ such that 
\[
 T=X\Delta_{1}^{-\frac{m}{2}}+ R
\]
where $R$ has analytic order less than or equal to $\ell$.

A finite linear combinations
of basic pseudodifferential operator of order $\leq k$ is called a
pseudodifferential operator of order $\leq k$.
\end{dfn}

We denote the set of pseudodifferential operators of order $\leq 0$ by
$\Psi_{0}(\mathcal{D},\Delta)$. It is proved in \cite{Higson_localindex} that
the pseudodifferential operators of order $\leq 0$ is infact an algebra. We
need the following proposition due to Higson. Denote the derivation $T \mapsto
[\Delta^{\frac{1}{2}},T]$ by $\delta$.

\begin{ppsn}
\label{pseudo}
Let $(\mathcal{D},\Delta)$ be a differential pair. The derivation $\delta$
leaves the algebra $\Psi_{0}(\mathcal{D},\Delta)$ invariant.
\end{ppsn}

Let $(\mathcal{D},\Delta)$ be a differential pair. Assume that
$\Delta^{-\frac{r}{2}}$ is trace class for some $r>0$. We say that the
analytic dimension of $\mathcal{D},\Delta)$ is $p$ if 
\[
 p:=\inf\{q>0: \Delta^{\frac{-r}{2}} \text{~ is trace class for every~} r>q \}
\]
 
 us make the
following definition of the heat kernel expansion for a differential pair. 
\begin{dfn}
 Let $(\mathcal{D},\Delta)$ be a differential pair of analytic dimension $p$.
 We say that $(\mathcal{D},\Delta)$ has a heat kernel expansion
if for $X \in \mathcal{D}_{m}$, the function $t\mapsto
t^{p+m}Tr(Xe^{-t^{2}\Delta})$ has an asymptotic expansion near $0$. 
\end{dfn}

Now we show that if $(\mathcal{D}.\Delta)$ has the heat kernel expansion then
the algebra $\Psi_{0}(\mathcal{D},\Delta)$ has the weak heat kernel expansion.

\begin{ppsn}
\label{heat kernel for a differential pair}
 Let $(\mathcal{D},\Delta)$ be a differential pair of analytic dimension
$p$ having the heat kernel expansion. Denote the operator $\Delta^{\frac{1}{2}}$
by $|D|$. Then for every $b \in \Psi_{0}(\mathcal{D},\Delta)$, the function $t
\mapsto t^{p}Tr(be^{-t|D|})$ has an asymptotic power series expansion.
\end{ppsn}
\textit{Proof.} First observe that if $R:\clh_{\infty} \to \clh_{\infty}$ is an
operator of analytic order $<-p-n-1$  then $R|D|^{n+1}$ is trace class and
hence by Taylor's series \[
                                            Tr(Re^{-t|D|})=
\sum_{k=0}^{n}\frac{(-1)^{k}Tr(R|D|^{k})}{k!}t^{k}+ O(t^{n+1})
                                           \]
for $t$ near $0$. Thus it is enough to show the result when
$b=X\Delta_{1}^{-\frac{m}{2}}$. For an operator $T$ on $\clh_{\infty}$, let
$\zeta_{T}(s):= Tr(X|D|^{-s})$. Then $\zeta_{b}(s):=\zeta_{X}(s+m)$. As in 
Proposition \ref{implication} one can show that $\Gamma(s)\zeta_{X}(s)$ has
rapid decay along the vertical strips. Now \[
\Gamma(s)\zeta_{b}(s)=\frac{\Gamma(s)}{\Gamma(s+m)}\Gamma(s+m)\zeta_{X}(s+m)
                                         \]
 Hence
$\Gamma(s)\zeta_{b}(s)$ has rapid decay along the vertical strips. But
$\Gamma(s)\zeta_{b}(s)$ is the Mellin transform of $Tr(be^{-t|D|})$. Hence by
Proposition \ref{IMTT}, it follows that $t^{p}Tr(be^{-t|D|})$ has an asymptotic
power series expansion. This completes the proof. \hfill $\Box$
      
We make use of the following proposition to prove that spectral triple
associated to the NC torus and that of a spin manifold posses the weak heat
kernel expansion property.

\begin{ppsn}
\label{diff implies heat kernel}
 Let $(\mathcal{A},\clh,D)$ be a finitely summable spectral triple and
$\Delta:=D^{2}$. Suppose that there exists an algebra of operators
$D:=\bigcup_{p \geq 0}\mathcal{D}_{p}$ such that $(D,\Delta)$ is a differential
pair of analytic dimension $p$. Assume that $(\mathcal{D},\Delta)$ satisfies the
following
\begin{enumerate}
 \item The algebra $\mathcal{D}_{0}$ contains $\mathcal{A}$ and
$[D,\mathcal{A}]$.
 \item The differential pair $(\mathcal{D},\Delta)$ has the heat kernel
expansion property.
\item The operator $D \in \mathcal{D}_{1}$.
\end{enumerate}
Then the spectral triple $(\mathcal{A},\clh,D)$ has the weak heat kernel
asymptotic expansion property.
\end{ppsn}
\textit{Proof.} Without loss of generality, we can assume that $D$ is
invertible. We let $\mathcal{B}$ be the algebra of pseudodifferential operators
of order $0$ associated to $(\mathcal{D},\Delta)$ . Now Proposition \ref{pseudo}
together with the fact that $\mathcal{D}_{0} \subset \mathcal{B}$ shows that
$\mathcal{B}$ contains $\mathcal{A}$ and $[D,\mathcal{A}]$ and is invariant
under $\delta:=[|D|,.]$. Since $D \in \mathcal{D}_{1}$, it follows that
$F:=D\Delta^{\frac{-1}{2}} \in \mathcal{B}$. Now $(4)$ of Definition \ref{heat
kernel} follows from Proposition \ref{heat kernel for a differential pair}. This
completes the proof. \hfill $\Box$

\subsection{Examples}
Now we discuss some examples of spectral triples which satisfy the weak heat
kernel asymptotic expansion. We start with the classical example.

 Let $M$ be a Reimannian spin manifold and $S\to M$ be a spinor bundle. We
denote the Hilbert space of square integrable sections on $L^{2}(M,S)$ by
$\clh$. We represent $C^{\infty}(M)$ on $\clh$ by multiplication operators. Let
$D$ be the Dirac
operator associated with Levi-Civita connection. Then 
the triple $(C^{\infty}(M),\clh,D)$ is a spectral triple. Then the operator
$D^{2}$ is then a generalised Laplacian ( Refer \cite{Berline}). Let
$\mathcal{D}$ denote the usual algebra of differential operators on $S$. Then
$(\mathcal{D},\Delta)$ is a differential pair. Moreover Proposition 2.4.6 in
\cite{Berline} implies that $(\mathcal{D},\Delta)$ has the heat kernel
expansion. Also $D \in \mathcal{D}_{1}$. Now Proposition \ref{diff implies heat
kernel} implies that the spectral triple $(C^{\infty}(M), \clh, D)$ has the
weak heat kernel asymptotic expansion.

\subsubsection{The spectral triple associated to the NC torus}
 Let us recall the definition of the Noncommutative torus which we abbreviate
as NC torus. Throughout we assume that $\theta \in [0,2\pi)$.
\begin{dfn}
 The $C^{*}$ algebra $A_{\theta}$ is defined as the universal $C^{*}$ algebra
generated by two unitaries $u$ and $v$ such that $uv=e^{i\theta}vu$
\end{dfn}
Define the operators $U$ and $V$ on $\ell^{2}(\mathbb{Z}^{2})$ as follows:
\begin{align*}
 Ue_{m,n}&:=e_{m+1,n}\\
Ve_{m,n}&:=e^{-in\theta}e_{m,n+1}
\end{align*}
where $\{e_{m,n}\}$ denotes the standard orthonormal basis on
$\ell^{2}(\mathbb{Z}^{2})$.
Then it is well known that $u\to U$ and $v \to V$ gives a faithful
representation of the $C^{*}$ algebra $A_{\theta}$. 

Consider the positive selfadjoint operator $\Delta$ on
$\clh:=\ell^{2}(\mathbb{Z}^{2})$ defined on the orthonormal basis $\{e_{m,n}\}$
by $\Delta(e_{m,n})= (m^{2}+n^{2})e_{m,n}$. For a polynomial $P=p(m,n)$, define
the operator $T_{P}$ on $\clh_{\infty}$ by $T_{P}(e_{m,n}):=p(m,n)e_{m,n}$.
The group $\mathbb{Z}^{2}$ acts on the algebra of polynomials as follows. For
$x:=(a,b) \in \mathbb{Z}^{2}$ and $P:=p(m,n)$, define
$x.P:=p(m-a,n-b)$. We denote $(1,0)$ by $e_{1}$ and $(0,1)$ by $e_{2}$.

Note that if $P$ is a polynomial of degree $\leq k$, then
$T_{P}\Delta^{-\frac{k}{2}}$ is bounded on $Ker(\Delta)^{\perp}$. Thus it
follows that if $P$ is a polynomial of degree $\leq k$ then $T_{P}$ has
analytic order $\leq k$.

Also note that \[
\Delta_{1}^{\frac{k}{2}}U\Delta_{1}^{\frac{-k}{2}}e_{m,n}:=
\frac{((m+1)^{2}+n^{2})^{\frac{k}{2}}}{(m^{2}+n^{2})^{\frac{k}{2}}}e_{m+1,n}
\text{~if~} (m,n)\neq 0
                \]
Thus it follows that $U$ is of analytic order $\leq 0$. Similary one can show
that $V$ is of analytic order $\leq 0$. Now note the following commuatation
relationship
\begin{align}
 UT_{P}:&=T_{e_{1}.P}U \\
 VT_{P}:&=T_{e_{1}.P}V
\end{align}
Thus it follows that $[\Delta,U^{\alpha}V^{\beta}]=T_{Q}U^{\alpha}V^{\beta}$
for some degree 1 polynomial $Q$. 

Let us define $\mathcal{D}_{p}:= span\{T_{P_{\alpha,\beta}}U^{\alpha}V^{\beta}:
deg(P_{\alpha,\beta}) \leq k \}$ and let
$\mathcal{D}:=\bigcup_{p}\mathcal{D}_{p}$. The above observations can be
rephrased into the following proposition.
\begin{ppsn}
 \label{diff pair for NC torus}
The pair $(\mathcal{D},\Delta)$ is a differential pair of analytic dimension
$2$.
\end{ppsn}

Now we show that the differential pair $(\mathcal{D},\Delta)$ has the heat
kernel expansion.
\begin{ppsn}
 \label{heat kernel for NC torus}
The differential pair $(\mathcal{D},\Delta)$ has the heat kernel expansion
property.
\end{ppsn}
\textit{Proof.} Let $X \in \mathcal{D}_{q}$ be given. It is enough to consider
the case when $X:=T_{P}U^{\alpha}V^{\beta}$. First note that
$Tr(Xe^{-t\Delta})=0$ unless $(\alpha,\beta)=0$. Now let $X:=T_{P}$. Again it
is enough to consider the case when $P$ is a monomial. Let
$P=p(m,n)=m^{k_{1}}n^{k_{2}}$. Now \[ Tr(T_{P}e^{-t\Delta})= (\sum_{m \in
\mathbb{Z}}m^{k_{1}}e^{-tm^{2}})(\sum_{n \in \mathbb{Z}}n^{k_{2}}e^{-tn^{2}})   
                                                                   \]
Now the asymptotic expansion follows from applying Proposition 2.4.6 in
\cite{Berline} to the standard Laplacian on the circle. This completes the
proof. \hfill $\Box$

Let $\mathcal{A}_{\theta}$ be the $*$ algebra generated by $U$ and $V$. We
consider the direct sum  representation of $A_{\theta}$ on
$\clh \oplus \clh$. Define $D:=\begin{bmatrix}
                                                           0 & T_{m-in} \\
                                                           T_{m+in} & 0
                                                          \end{bmatrix}$. Then
$D$ is selfadjoint on $\clh \oplus \clh$ and $D^{2}=\begin{bmatrix}
                                                     \Delta & 0 \\
                                                      0 & \Delta
                                                   \end{bmatrix}$.
 It is well  known that $(\mathcal{A}_{\theta}, \clh \oplus \clh, D)$ is a
$2+$ summable spectral  triple. 

\begin{ppsn}
 The spectral triple $(\mathcal{A}_{\theta}, \clh\oplus \clh, D)$ has the weak
heat kernel asymptotic expansion property.
\end{ppsn}
\textit{Proof.} Let  $(\mathcal{D},\Delta)$ be the differential pair considered
in Proposition \ref{diff pair for NC torus}. Then the amplification
$(\mathcal{D}^{'}:=M_{2}(\mathcal{D}),D^{2})$ is a differential pair. Note that
$D \in \mathcal{D}^{'}_{1}$. Clearly $\mathcal{A}_{\theta} \in
\mathcal{D}^{'}$. Note the commutation relations
\begin{align*}
 [T_{m \pm in},U]&= U \\
 [T_{m \pm in},V]&=\pm iV
\end{align*}

This implies that $[D,\mathcal{A}_{\theta}] \subset \mathcal{D}^{'}_{0}$. Since
$(\mathcal{D},\Delta)$ has the heat kernel expansion, it follows that the
differential pair $(M_{2}(\mathcal{D}),D^{2})$ also has the heat kernel
expansion. Now Proposition \ref{diff implies heat kernel} implies that the
spectral triple $(\mathcal{A}_{\theta},\clh \oplus \clh, D)$ has the weak heat
kernel expansion. This completes the proof. \hfill $\Box$

\subsubsection{The torus equivariant spectral triple on the odd dimensional
quantum spheres}
In this section we recall the spectral triple for the odd dimensional quantum
spheres given in \cite{PsPal1}.
We begin with some known facts about odd dimensional quantum spheres. 
Let $q\in(0,1]$.
The $C^*$-algebra $C(S_q^{2\ell+1})$ of the quantum
sphere $S_q^{2\ell+1}$
is the universal $C^*$-algebra generated by
elements
$z_1, z_2,\ldots, z_{\ell+1}$
satisfying the following relations (see~\cite{Hong-Sym1}):
\bean
z_i z_j & =& qz_j z_i,\qquad 1\leq j<i\leq \ell+1,\\
z_i^* z_j & =& q z_j z_i^* ,\qquad 1\leq i\neq j\leq \ell+1,\\
z_i z_i^* - z_i^* z_i +
(1-q^{2})\sum_{k>i} z_k z_k^* &=& 0,\qquad \hspace{2em}1\leq i\leq \ell+1,\\
\sum_{i=1}^{\ell+1} z_i z_i^* &=& 1.
\eean
We will denote by $\cla(S_q^{2\ell+1})$ the *-subalgebra of $A_\ell$
generated by the $z_j$'s. Note that for $\ell=0$, the $C^*$-algebra
$C(S_q^{2\ell+1})$ is the algebra of continuous functions
$C(\bbt)$ on the torus and for $\ell=1$, it is $C(SU_q(2))$.

There is a natural torus group $\mathbb{T}^{\ell +1}$ action $\tau$ on
$C(S_{q}^{2\ell+1})$ as follows. For $w=(w_{1},\ldots,w_{\ell +1})$, define
an automorphism $\tau_{w}$ by $\tau_{w}(z_{i})=w_{i}z_{i}$.  

Recall that $N$ is
the number operator on $\ell^{2}(\bbn)$ and $S$ is the left
shift on $\ell^{2}(\bbn)$. We also use the same notation $S$ for the left shift
on $\ell^{2}(\mathbb{Z})$. We let $\clh_{\ell}$ denote the Hilbert space
$\ell^{2}(\bbn^{\ell} \times Z)$.
Let $Y_{k,q}$ be the following operators on $\clh_\ell$:
\be\label{eq:ykq}
 Y_{k,q}=\begin{cases}
 \underbrace{q^N\otimes\ldots\otimes q^N}_{k-1 \mbox{
copies}}\otimes
      \sqrt{1-q^{2N}}S^*\otimes 
   \underbrace{I \otimes\cdots\otimes I}_{\ell+1-k \mbox{ copies}}, & \mbox{ if
} 1\leq k\leq \ell,\cr
   &\cr
    \underbrace{q^N\otimes\cdots\otimes q^N}_{\ell \mbox{ copies}}
       \otimes S^*, &  \mbox{ if } k=\ell+1.
         \end{cases}
\ee

Then $\pi_\ell:z_k\mapsto Y_{k,q}$ gives a faithful representation
of $C(S_q^{2\ell+1})$ on $\clh_\ell$ for $q\in(0,1)$ 
(see lemma~4.1 and remark~4.5, \cite{Hong-Sym1}).
We will denote the image $\pi_\ell(C(S_q^{2\ell+1}))$ by $A_\ell(q)$ or by just
$A_\ell$.

 Let  $\{e_{\gamma}: \gamma \in \Gamma_{\Sigma_\ell} \}$ 
be the standard orthonormal basis for $\clh_{\ell }$. We recall the following
theorem from \cite{PsPal1}.

\bthm[\cite{PsPal1}]  
Let $D_{\ell}$ be the operator $e_{\gamma} \to d(\gamma)e_{\gamma}$ on
$\clh_{\ell }$
where the $d_{\gamma}$'s are given by
\begin{displaymath} 
\begin{array}{lll}
d(\gamma)&=&\left\{\begin{array}{ll}
                            \gamma_{1}+\gamma_{2}+\cdots
\gamma_{\ell }+|\gamma_{\ell +1}| & \text{ if  }\gamma_{\ell +1} \geq 0 \\
                            -(\gamma_{1}+\gamma_{2}+\cdots
\gamma_{\ell }+|\gamma_{\ell +1}|) & \text{ if }   \gamma_{\ell +1} < 0
                       \end{array} \right.
\end{array}
\end{displaymath}
Then $(\cla(S_{q}^{2\ell+1}),\clh_{\ell },D_{\ell})$ is a non-trivial $(\ell+1)$
summable spectral triple. 
\ethm

But to deduce  that the spectral triple
$(\mathcal{A}(S_{q}^{2\ell+1}), \clh_{\ell}, D_{\ell})$ satisfies the weak heat
kernel asymptotic expansion, we need a topological version of
Definition \ref{heat kernel} and Proposition \ref{suspension}. We do this in
the next section.

\section{Smooth subalgebras and the weak heat kernel asymptotic expansion}
First we recall the definition of smooth subalgebras of $C^{*}$ algebras. For an
algebra $A$ (possibly non-unital), we denote the algebra obtained by adjoining a
unit to $A$ by $A^{+}$.
\begin{dfn}
 Let $A$ be a unital $C^{*}$ algebra. A dense unital $*$ subalgebra
$\mathcal{A}^{\infty}$ is called smooth in $A$ if 
\begin{enumerate}
 \item The algebra $\mathcal{A}^{\infty}$ is a Fr\'echet $*$ algebra.
 \item The unital inclusion $\mathcal{A}^{\infty} \subset A$ is continuous.
 \item The algebra $\mathcal{A}^{\infty}$ is spectrally invariant in $A$ i.e. if
an element $a \in \mathcal{A}^{\infty}$ is invertible in $A$ then $a^{-1} \in
\mathcal{A}^{\infty}$.
\end{enumerate}
Suppose $A$ is a non-unital $C^{*}$ algebra. A dense Fr\'echet $*$ subalgebra
$\mathcal{A}^{\infty}$ is said to be smooth in $A$ if
$(\mathcal{A}^{\infty})^{+}$ is smooth in $A^{+}$.
\end{dfn}
 We also assume that our smooth subalgebras satisfy the condition that if
$\mathcal{A}^{\infty}\subset A$ is smooth  then
$\mathcal{A}^{\infty}\hat{\otimes}_{\pi}
\mathcal{S}(\ell^{2}(\mathbb{N}^{k}))\subset A \otimes
\mathcal{K}(\ell^{2}(\mathbb{N}^{k}))$ is smooth. 

Let $A$ be a unital $C^{*}$ algebra and $\mathcal{A}^{\infty}$ be a smooth
unital $*$ subalgebra of $A$. Assume that the topology on $\mathcal{A}^{\infty}$
is given by the countable family of seminorms $(\| \cdot \|_{p})$. Let us denote
the operator $1\otimes S$ by $\alpha$. Define the smooth quantum double
suspension of
$\mathcal{A}^{\infty}$ as follows 
\bea
\Sigma^{2}(\mathcal{A}^{\infty})&:= & \left\{ \sum_{j,k \in
\mathbb{N}}\alpha^{*j}(a_{jk}\otimes p)\alpha^{k} + \sum_{k\geq
0}\lambda_{k}\alpha^{k}+
\sum_{k >0}\lambda_{-k}\alpha^{*k}~:~a_{jk} \in
\mathcal{A}^{\infty},\right.\nonumber\\
&& \hspace{2em}\left.\sum_{j,k}(1+j+k)^{n} \| a_{jk}
\|_{p}<\infty,~(\lambda_{k})
\text{ is rapidly decreasing} \right\}.\label{eq:Aellinfty}
\eea
Now let us topologize $\Sigma^{2}(\mathcal{A}^{\infty})$ by defining a seminorm
$\|~\|_{n,p}$ for every $n,p \geq 0$. For an element 
\begin{displaymath}
a:=\sum_{j,k \in \mathbb{N}}\alpha^{*j}(a_{jk}\otimes p)\alpha^{k} + \sum_{k\geq
0}\lambda_{k}\alpha^{k}+
\sum_{k >0}\lambda_{-k}\alpha^{*k} 
\end{displaymath}
 in $\Sigma^{2}(\mathcal{A}^{\infty})$ we define $\|a\|_{n,p}$ by 
\begin{displaymath}
 \|a\|_{n,p}:= \sum_{j,k \in \mathbb{N}}(1+|j|+|k|)^{n}\|a_{jk}\|_{p}+\sum_{k
\in \mathbb{Z}}(1+|k|)^{n}|\lambda_{k}|
\end{displaymath}
It is easily verifiable that 
\begin{enumerate}
 \item The subspace $\Sigma^{2}(\mathcal{A}^{\infty})$ is a dense $*$ subalgebra
of $\Sigma^{2}(A)$.
 \item The topology on $\Sigma^{2}(\mathcal{A}^{\infty})$ induced by the
seminorms $(\|~\|_{n,p})$ makes $\Sigma^{2}(\mathcal{A}^{\infty})$ a Fr\'echet
$*$ algebra.
 \item The unital inclusion $\Sigma^{2}(\mathcal{A}^{\infty}) \subset 
\Sigma^{2}(A)$ is continuous.
\end{enumerate}
The next proposition proves that the Fr\'echet algebra
$\Sigma^{2}(\mathcal{A}^{\infty})$ is infact smooth in $\Sigma^{2}(A)$.
\begin{ppsn}
 Let $A$ be a unital $C^{*}$ algebra and let $\mathcal{A}^\infty \subset A$ be a
unital smooth subalgebra such that $\mathcal{A}^{\infty}\hat{\otimes}_{\pi}
\mathcal{S}(\ell^{2}(\mathbb{N}^{k}))\subset A \otimes
\mathcal{K}(\ell^{2}(\mathbb{N}^{k}))$ is smooth for every $k \in \mathbb{N}$.
Then the algebra $\Sigma^{2}(\mathcal{A}^{\infty}) \hat{\otimes}_{\pi}
\mathcal{S}(\ell^{2}(\mathbb{N}^{k}))$ is smooth in $\Sigma^{2}(A)\otimes
\mathcal{K}(\ell^{2}(\mathbb{N}^{k}))$ for every $k\geq 0$.
Z\end{ppsn}
\textit{Proof.} Let us denote the restriction of $1\otimes \sigma$ to
$\Sigma^{2}(A)$ by $\rho$. Recall the $\sigma:\mathscr{T} \to C(\mathbb{T})$ is
the symbol map sending $S$ to the generating unitary.
Then one has the following exact sequence at the $C^{*}$ algebra level
\begin{displaymath}
 0 \to A \otimes \mathcal{K}(\ell^{2}(\mathbb{N})) \to \Sigma^{2}(A)
\stackrel{\rho} \to C(\mathbb{T}) \to 0
\end{displaymath}
At the subalgebra level one has the following ``sub'' exact sequence
\begin{displaymath}
 0 \to \mathcal{A}^{\infty} \hat{\otimes}_{\pi}
\mathcal{S}(\ell^{2}(\mathbb{N})) \to \Sigma^{2}(\mathcal{A}^{\infty})
\stackrel{\rho} \to C^{\infty}(\mathbb{T}) \to 0
\end{displaymath}
Since $\mathcal{A}^{\infty} \hat{\otimes}_{\pi}\mathcal{S}(\ell^{2}(\mathbb{N}))
\subset A \otimes \mathcal{K}(\ell^{2}(\mathbb{N}))$ and $C^{\infty}(\mathbb{T})
\subset C(\mathbb{T})$ are smooth, it follows from theorem~3.2, part ~2,
\cite{Sch-1993} that $\Sigma^{2}(\mathcal{A}^{\infty})$ is smooth in
$\Sigma^{2}(A)$. One can prove  that $\Sigma^{2}(\mathcal{A}^{\infty})
\hat{\otimes}_{\pi} \mathcal{S}(\ell^{2}(\mathbb{N}^{k}))$ is smooth in
$\Sigma^{2}(A)\otimes \mathcal{K}(\ell^{2}(\mathbb{N}^{k}))$ for every $k>0$
along the same lines first by tensoring the $C^{*}$ algebra exact sequence by
$\mathcal{K}(\ell^{2}(\mathbb{N}^{k}))$ and then by tensoring the Fr\'echet
algebra exact sequence by $\mathcal{S}(\ell^{2}(\mathbb{N}^{k}))$ and appealing
to Theorem~3.2,part~2 of \cite{Sch-1993}. This completes the proof. \hfill
$\Box$
\subsection{ The topological weak heat kernel expansion }
We need the following version of the weak heat kernel expansion.
\begin{dfn}
\label{topological heat kernel}
Let $(\mathcal{A}^{\infty},\mathcal{H},D)$ be a $p+$ summable spectral
triple for a $C^{*}$ algebra A where $\mathcal{A}^{\infty}$ is smooth in $A$. We
say that the spectral triple $(\mathcal{A}^{\infty},\mathcal{H},D)$ has the
topological weak heat kernel asymptotic expansion property if there exists a $*$
subalgebra $\mathcal{B}^{\infty} \subset B(\mathcal{H})$ such that 
\begin{itemize}
\item[(1)] The algebra $\mathcal{B}^{\infty}$ has a Fr\'echet space structure
and endowed with it it is a Fr\'echet $*$ algebra.
\item[(2)] The algebra $\clb^{\infty}$ contains
$\mathcal{A}^{\infty}$.
\item[(3)] The inclusion $\mathcal{B}^{\infty}\subset B(\clh)$ is continuous.
\item[(4)] The unbounded derivations $\delta:=[|D|,.]$ leaves
$\mathcal{B}^{\infty}$ invariant and is continuous. Also the unbounded
derivation $d:=[D,.]$ maps $\mathcal{A}^{\infty}$ into $\mathcal{B}^{\infty}$
in a continuous fashion.
\item[(5)] The left multiplication by the operator $F:=sign(D)$ denoted $L_{F}$ 
leaves $\mathcal{B}^{\infty}$ invariant and is continuous. 
\item[(6)] The function $\tau_{p}:(0,\infty) \times \mathcal{B}^{\infty} \to
\mathbb{C}$ defined by $\tau_{p}(t,b)=t^{p}Tr(be^{-t|D|})$ has a uniform
asymptotic power series expansion.
\end{itemize}
\end{dfn}

We need the following analog of Proposition \ref{tensoring by smooth compacts}
and Proposition \ref{suspension}. First we need the following two lemmas.

\begin{lmma}
 \label{extension of asymptotic expansion}
Let $E$ be a Fr\'echet space and $F \subset E$ be a dense subspace. Let
$\phi:(0,\infty)\times E \to \mathbb{C}$ be a continuous function which is
linear in the second variable. Suppose that $\phi:(0,\infty)\times F \to
\mathbb{C}$ has a uniform asymptotic power series expansion then
$\phi:(0,\infty)\times E \to \mathbb{C}$ has a uniform asymptotic power series
expansion.
\end{lmma}
\textit{Proof.} Suppose that $\phi(t,f) \sim \sum_{r=0}^{\infty}a_{r}(f)t^{r}$.
Then $a_{r}:F\to\mathbb{C}$ is linear and is continuous for every $r\in
\mathbb{N}$. Since $F$ is dense in $E$, for every $r \in \mathbb{N}$, the
function $a_{r}$ admits a continuous extension to the whole of $E$ which we
still denote it by $a_{r}$. Now fix $N \in \mathbb{N}$. Then there exists a
neighbourhood $U$ of $E$ containing $0$ and $\epsilon ,M >0$ such that 
\begin{equation}
\label{dense}
 |\phi(t,f)-\sum_{r=0}^{N}a_{r}(f)t^{r}| \leq M t^{N+1} \textrm{ ~for~}
0<t<\epsilon, f \in U \cap F  
\end{equation}
Since $\phi(t,.)$ and $a_{r}(.)$ are continuous and as $F$ is dense in $E$,
Equation \ref{dense} continues to hold for every $f \in U$. This completes the
proof. \hfill $\Box$
\begin{lmma}
\label{tensor product}
 Let $E_{1},E_{2}$ be Fr\'echet spaces and let $F_{i}:(0,\infty)\times E_{i} \to
\mathbb{C}$ be continuous and linear in the second variable for $i=1,2$.
Consider the function $F:(0,\infty)\times E_{1}\hat{\otimes}_{\pi} E_{2} \to
\mathbb{C}$ be defined by $F(t,e_{1}\otimes e_{2})=F_{1}(t,e_{1})F(t,e_{2})$.
Assume that $F$ is continuous. If $F_{1}$ and $F_{2}$ have uniform asymptotic
expansions then $F$ has a uniform asymptotic power series expansion.
\end{lmma}
\textit{Proof.} By Lemma \ref{extension of asymptotic expansion}, it is enough
to show that $F:(0,\infty)\times E_{1}\otimes_{alg}E_{2} \to \mathbb{C}$ has a
uniform asymptotic power series expansion. Let $\theta:E_{1} \times E_{2} \to
E_{1}\otimes_{alg} E_{2}$ be defined by $\theta(e_{1},e_{2})=e_{1}\otimes
e_{2}$. Consider the map $G:(0,\infty)\times E_{1} \times E_{2} \to \mathbb{C}$ 
defined by $G(t,e_{1},e_{2}):=F(t,\theta((e_{1},e_{2})))$. By Lemma
\ref{product}, it follows that $G$ has a uniform asymptotic power series
expansion say 
\begin{displaymath}
G(t,e) \sim \sum_{r=0}^{\infty}a_{r}(e)t^{r}
\end{displaymath}
 The maps $a_{r}:E_{1}\times E_{2} \to \mathbb{C}$ are continuous bilinear. We
let $\tilde{a}_{r}:E_{1}\hat{\otimes}_{\pi}E_{2} \to \mathbb{C}$ be the linear
maps such that $\tilde{a}_{r}\circ \theta:=a_{r}$. Let $N \in \mathbb{N}$ be
given. Then there exists $\epsilon,M>0$ and open sets $U_{1},U_{2}$ containing
$0$ in $E_{1},E_{2}$ such that 
\begin{equation}
\label{estimate}
 |G(t,e)-\sum_{r=0}^{N}a_{r}(e)t^{r}| \leq M t^{N+1} \textrm{ ~for~}
0<t<\epsilon, e \in U_{1}\times U_{2}   
\end{equation}
Without loss of generality, we can assume that $U_{i}:=\{x \in
E_{i}:~p_{i}(x)<1\}$ for a seminorm $p_{i}$ of $E_{i}$. Now Equation
\eqref{estimate} implies that 
\begin{equation}
 |F(t,\theta(e))-\sum_{r=0}^{N}\tilde{a}_{r}(\theta(e))t^{r}| \leq M t^{N+1}
\textrm{ ~for~} 0<t<\epsilon, e \in U_{1}\times U_{2} 
\end{equation}
 Hence for $t\in(0,\epsilon)$ and $x \in \theta(U_{1}\times U_{2})$,
\begin{equation} 
\label{estimate1}
|F(t,x)-\sum_{r=0}^{N}\tilde{a}_{r}(x)t^{r}| \leq M t^{N+1} 
\end{equation}
Since $\tilde{a}_{r}$ is linear and $F$ is linear in the second variable, it
follows that Equation \ref{estimate1} continues to hold for $x$ in the convex
hull of $\theta(U_{1}\times U_{2})$ which is nothing but the unit ball
determined by the seminorm $p_{1}\otimes p_{2}$ in $E_{1}\otimes_{alg}E_{2}$.
This completes the proof. \hfill $\Box$

In the next proposition, we consider the stability of the weak heat kernel
asymptotic
expansion property for tensoring by smooth compacts.
\begin{ppsn}
\label{topological tensoring by smooth compacts}
 Let $(\mathcal{A}^{\infty},\mathcal{H},D)$ be a spectral triple where the
algebra $\mathcal{A}^{\infty}$ is a smooth subalgebra of $C^{*}$ algebra.
Assume that $(\mathcal{A}^{\infty},\mathcal{H},D)$ has the topological weak
heat kernel expansion property with dimension $p$. Then the spectral triple
$(\mathcal{A}^{\infty}\hat{\otimes}_{\pi}
\mathcal{S}(\ell^{2}(\mathbb{N})),\mathcal{H}\otimes
\ell^{2}(\mathbb{N}),D_{0}:=(F\otimes 1)(|D|\otimes 1+1\otimes N))$ also has the
weak heat kernel asymptotic expansion property with dimension $p$ where
$F:=sign(D)$.
\end{ppsn}
\textit{Proof.} Let $\mathcal{B}^{\infty} \subset B(\mathcal{H})$ be a $*$
subalgebra for which $(1)-(6)$ of Definition \ref{topological heat kernel} are
satisfied. We
denote
$\mathcal{B}^{\infty}\hat{\otimes}_{\pi}\mathcal{S}(\ell^{2}(\mathbb{N}))$ by
$\mathcal{B}_{0}^{\infty}$. We show that $\mathcal{B}_{0}^{\infty}$ satisfies
$(1)-(6)$ of Definition \ref{topological heat kernel}. First note that the
natural
representation of $\mathcal{B}_{0}^{\infty}$ in $\mathcal{H}\otimes
\ell^{2}(\mathbb{N})$ is injective. Thus $(3)$ is clear. Also $(1)$ and $(2)$
are obvious. Now let us now prove $(4)$.

 We denote the unbounded derivation $[|D_{0}|,.]$,$[|D|,.]$ and $[N,.]$ by
$\delta_{D_{0}},\delta_{D}$ and $\delta_{N}$ respectively. By assumption
$\delta_{D}$ leaves $\mathcal{B}$ invariant and is continuous. It is also easy
to see that $\delta_{N}$ leaves $\mathcal{S}(\ell^{2}(\mathbb{N}))$ invariant
and is continuous. Let $\delta^{'}:=\delta_{D}\otimes 1+1\otimes \delta_{N}$.
Then $\delta^{'}:\mathcal{B}_{0}^{\infty} \to \mathcal{B}_{0}^{\infty}$ is
continuous. Clearly $\mathcal{B}^{\infty}\otimes_{alg}
\mathcal{S}(\ell^{2}(\mathbb{N}))$ is contained in the domain of $\delta$ and
$\delta=\delta^{'}$ on $\mathcal{B}^{\infty}\otimes_{alg}
\mathcal{S}(\ell^{2}(\mathbb{N}))$. Now let $a \in \mathcal{B}_{0}^{\infty}$ be
given. Then there exists a sequence $(a_{n})$ in
$\clb^{\infty}\otimes_{\pi}\cls(\ell^{2}(\bbn))$ such that $(a_{n})$ converges
to $a$ in $\clb_{0}^{\infty}$. Since $\delta^{'}$ is continuous on
$\clb_{0}^{\infty}$ and the inclusion $\clb_{0}^{\infty} \subset B(\clh)$ is
continuous, it follows that $\delta_{D_{0}}(a_{n})=\delta^{'}(a_{n})$ converges
to $\delta^{'}(a)$. As $\delta_{D_{0}}$ is a closed derivation, it follows that
$a \in Dom(\delta_{D_{0}})$ and $\delta_{D_{0}}(a)=\delta^{'}(a)$. Hence we have
shown that $\delta_{D_{0}}$ leaves $\clb_{0}^{\infty}$ invariant and is
continuous. Similarly one can show that the unbounded derivation
$d_{0}:=[D_{0},.]$
maps $\cla\hat{\otimes}_{\pi}S(\ell^{2}(\bbn))$ into $\clb_{0}^{\infty}$
invariant in a continuous manner. 

As $F_{0}:=sign(D_{0})=F\otimes 1$, $(5)$ is clear.
 Consider the function $\tau_{p}:(0,\infty)\times \clb_{0}^{\infty} \to \bbc$
defined by $\tau_{p}(t,b):=t^{p}Tr(be^{-t|D_{0}|})$. Then $\tau_{p}(t,b\otimes
k)=\tau_{p}(t,b)\tau_{0}(t,k)$. Hence by Lemma \ref{tensor product}, it follows
that $\tau_{p}$ has a uniform asymptotic power series expansion. This completes
the proof. \hfill $\Box$

Now we consider the stability of the weak heat kernel asymptotic expansion under
the
double suspension. 

\begin{ppsn}
\label{topological suspension}
 Let $(\mathcal{A}^{\infty},\mathcal{H},D)$ be a spectral triple with the
topological weak heat
kernel asymptotic expansion property of dimension $p$. Assume that the algebra
$\mathcal{A}^{\infty}$ is unital and the representation on $\mathcal{H}$ is
unital. Then the spectral triple $
(\Sigma^{2}(\mathcal{A}^{\infty}),\mathcal{H}\otimes
\ell^{2}(\mathbb{N}),\Sigma^{2}(D))$ also has
the topological weak heat kernel asymptotic expansion property with dimension
$p+1$.
\end{ppsn}
\textit{Proof.} We denote the operator $\Sigma^{2}(D)$ by $D_{0}$. Let
$\mathcal{B}^{\infty}$ be $*$ subalgebra of $B(\mathcal{H})$
for which $(1)-(6)$ of Definition \ref{topological heat kernel} are satisfied.
For
$f=\sum_{n \in \mathbb{Z}}\lambda_{n}z^{n} \in C(\mathbb{T})$, we let
$\sigma(f):= \sum_{n\geq 0}\lambda_{n}S^{n}+\sum_{n>0}\lambda_{-n}S^{*n}$.  We
denote the projection $\frac{1+F}{2}$ by $P$. We assume here that $P \neq \pm 1$
as the case $P=\pm 1$ is similar. We let $\mathcal{B}_{0}^{\infty}$ to denote
the algebra $\mathcal{B}^{\infty}\hat{\otimes}_{\pi}\cls(\ell^{2}(\bbn))$ as in
Proposition \ref{topological tensoring by smooth compacts}. As in Proposition
\ref{topological tensoring
by smooth compacts}, we let $\delta_{D_{0}},\delta_{D},\delta_{N}$ to denote the
unbounded derivations $[|D_{0}|,.],[|D|,.]$ and $[N,.]$ respectively. Define
\begin{displaymath}
 \mathcal{\tilde{B}}^{\infty}:=\{b+P\otimes \sigma(f)+(1-P)\otimes \sigma(g): ~b
\in \mathcal{B}_{0}^{\infty},f,g \in C^{\infty}(\mathbb{T})\}
\end{displaymath}
Then $\mathcal{\tilde{B}}^{\infty}$ is isomorphic to the direct sum
$\mathcal{B}_{0}^{\infty}\oplus C^{\infty}(\mathbb{T})\oplus
C^{\infty}(\mathbb{T})$. We give $\mathcal{\tilde{B}}^{\infty}$ the Fr\'echet
space structure coming from this decomposition. It is easy to see that
$\mathcal{\tilde{B}}^{\infty}$ is a Fr\'echet $*$ subalgebra of
$B(\mathcal{H}\otimes \ell^{2}(\bbn))$. Clearly $(\pi \otimes
1)(\Sigma^{2}(\cla^{\infty})) \subset \mathcal{\tilde{B}}^{\infty}$. Thus we
have shown that $(1)$ and $(2)$ of Definition \ref{topological heat kernel} are
satisfied. Since $\clb_{0}^{\infty}$ is represented injectively on $\clh
\otimes \ell^{2}(\bbn)$, it follows that $\tilde{\clb}$ satisfies $(3)$.

We have already shown in Proposition \ref{topological tensoring by smooth
compacts} that
$\mathcal{B}_{0}^{\infty}$ is closed under $\delta_{D_{0}}$ and is continuous.
Also we have shown that  $d_{0}:=[D_{0},.]$ maps $\cla
\hat{\otimes}_{\pi}S(\ell^{2}(\bbn))$ into $\clb_{0}^{\infty}$
continuously. Now
note that 
\begin{align*}
\delta_{D_{0}}(P\otimes \sigma(f))&=P\otimes \sigma(if^{'})\\
\delta_{D_{0}}((1-P)\otimes \sigma(g))&=(1-P)\otimes \sigma(ig^{'}) \\
[D_{0},P\otimes \sigma(f)]&= P \otimes \sigma(if^{'}) \\
[D_{0},(1-P)\otimes \sigma(g)]&=-(1-P)\otimes \sigma(ig^{'})
\end{align*}
Thus it follows that $\delta_{D_{0}}$ leaves 
$\mathcal{\tilde{B}}^{\infty}$ invariant and is continuous. Also, it follows
that $d_{0}:=[D_{0},.]$ maps $\Sigma^{2}(\cla^{\infty})$ into $\tilde{\clb}$ in
a continuous manner.

Since $F_{0}:=sign(D_{0})=F\otimes 1$, it follows from definition that $F_{0}
\in \mathcal{\tilde{B}}^{\infty}$. Now we show that
$\mathcal{\tilde{B}}^{\infty}$ satisfies $(6)$.

We have already shown in Proposition \ref{topological tensoring by smooth
compacts} that the
function $\tau_{p}:(0,\infty)\otimes \clb_{0}^{\infty} \to \bbc$ defined by
$\tau_{p}(t,b):=t^{p}Tr(be^{-t|D_{0}|})$ has a uniform asymptotic power series
expansion. Hence $\tau_{p+1}$ restricted to $\mathcal{B}_{0}^{\infty}$ has a
uniform asymptorm asymptotic power series expansion. Now note that 
\begin{eqnarray}
\label{tasy}
 \tau_{p+1}(P\otimes \sigma(f))&=& (\int f(\theta)d\theta)
t^{p}Tr(Pe^{-t|D|})tTr(e^{-tN})\\
\label{tasy1}
 \tau_{p+1}((1-P)\otimes \sigma(g))&=&(\int
g(\theta)d\theta)t^{p}Tr((1-P)e^{-t|D|})tTr(e^{-tN})) 
\end{eqnarray}
Now recall that we have assumed that $\mathcal{A}^{\infty}$ is unital and hence
$P \in \mathcal{B}^{\infty}$. Hence $t^{p}Tr(xe^{-t|D|})$ has an asymptotic
power series expansion for $x\in \{P,1-P\}$. Also $tTr(e^{-tN})$ has an
asymptotic power series expansion. Now Equation \ref{tasy} and Equation
\ref{tasy1}, together with the earlier observation that $\tau_{p+1}$ restricted
to $\mathcal{B}_{0}^{\infty}$ has a uniform asymptotic power series expansion,
imply that the function $\tau_{p+1}:(0,\infty) \times
\mathcal{\tilde{B}}^{\infty} \to \bbc$ has a uniform asymptotic power series
expansion. This completes the proof. \hfill $\Box$

\begin{rmrk}
 If we start with the canonical spectral triple
$(C^{\infty}(\mathbb{T}),L^{2}(\mathbb{T}),\frac{1}{i}\frac{d}{d\theta})$ on the
circle  and apply the double suspension recursively one obtains the torus
equivariant spectral triple for the odd dimensional quantum spheres constructed
in \cite{PsPal1}. Now proposition 3.4 in \cite{SunPal} along with Theorem
\ref{topological suspension} implies that the torus equivariant spectral triple
on $S_{q}^{2\ell+1}$ satisfies the weak heat kernel expansion. Also Theorem
\ref{topological suspension}, along with Theorem
\ref{dimension spectrum}, gives a proof of Proposition 3.9 in \cite{SunPal}. 
\end{rmrk}
\subsection{ The equivariant spectral triple on $S_{q}^{2\ell+1}$}
 In this section, we show that the equivariant spectral triple on
$S_{q}^{2\ell+1}$ constructed in \cite{PsPal2} has the topological weak heat
kernel asymptotic expansion. First let us recall that the odd dimensional
quantum spheres can be realised as the quantum homogeneous space. Throughout we
assume $q \in (0,1)$. The $C^{*}$ algebra of the quantum group $SU_{q}(n)$
denoted by $C(SU_{q}(n))$ is  defined as the universal $C^{*}$ algebra generated
by $\{u_{ij}: 1 \leq i,j \leq \ell\}$ satisfying the following conditions 
\begin{equation}
 \sum_{k=1}^{n}u_{ik}u_{jk}^{*}=\delta_{ij}~~,~~ \sum_{k=1}^{n}u_{ki}^{*}u_{kj}
=\delta_{ij}
\end{equation}
\begin{equation}
\sum_{i_{1}=1}^{n}\sum_{i_{2}=1}^{n} \cdots \sum_{i_{n}=1}^{n}
E_{i_{1}i_{2}\cdots i_{n}}u_{j_{1}i_{1}}\cdots u_{j_{n}i_{n}} =
E_{j_{1}j_{2}\cdots j_{n}}
\end{equation}
where
\begin{displaymath}
\begin{array}{lll}
E_{i_{1}i_{2}\cdots i_{n}}&:=& \left\{\begin{array}{lll}
                                 0 &if& i_{1},i_{2},\cdots i_{n} \text{~are not
distinct} \\
                                 (-q)^{\ell(i_{1},i_{2},\cdots,i_{n})}  
                                 \end{array} \right. 
\end{array}
\end{displaymath}
where for a permutation $\sigma$ on $\{1,2,\cdots,n \}$, $\ell(\sigma)$ denotes
its length. The $C^{*}$ algebra has the quantum group structure with the
comultiplication being defined by 
\begin{equation*}
 \Delta(u_{ij}):= \sum_{k}u_{ik}\otimes u_{kj}
\end{equation*}

 Call the generators of $SU_{q}(n-1)$ as $v_{ij}$. The
map $\phi:C(SU_{q}(n)) \to C(SU_{q}(n-1))$ defined by
 \begin{eqnarray}
\label{subgroup}
 \phi(u_{ij}) &:=& \left\{\begin{array}{ll}
                        v_{i-1,j-1}&\text{if}~ 2 \leq i,j \leq n \\
                        \delta_{ij} & \text{otherwise}
                        \end{array} \right.
\end{eqnarray}
is a surjective unital $C^{*}$ algebra homomorphism such that $\Delta \circ \phi
= (\phi \otimes \phi)\Delta$. In this way the quantum group  $SU_{q}(n-1)$ is a
subgroup of the quantum group $SU_{q}(n)$. The $C^{*}$ algebra of the quotient
$SU_{q}(n)/SU_{q}(n-1)$ is defined as 
\begin{displaymath}
 C(SU_{q}(n)/SU_{q}(n-1)):=\{a \in C(SU_{q}(n)):~(\phi \otimes 1)\Delta(a)=
1\otimes a\}.
\end{displaymath}
  Also the $C^{*}$ algebra $C(SU_{q}(n)/SU_{q}(n-1))$ is 
generated by   $\{u_{1j}: ~1 \leq j \leq n\}$. Moreover the map 
$\psi:C(S_{q}^{2n-1}) \to C(SU_{q}(n)/SU_{q}(n-1))$ defined by
$\psi(z_{i}):=q^{-i+1}u_{1i}$ is an isomorphism.

Let $h$ be the Haar state on the quantum group $C(SU_{q}(\ell+1))$ and let
$L^{2}(SU_{q}(\ell+1))$ be the corresponding $GNS$ space. We denote the closure
of
$C(S_{q}^{2\ell+1})$ in $L^{2}(SU_{q}(\ell+1))$ by $L^{2}(S_{q}^{2\ell+1})$.
Then
$L^{2}(S_{q}^{2\ell+1})$ is invariant under the regular representation of
$SU_{q}(\ell+1)$. Thus we get a covariant representation for the dynamical
system
$(C(S_{q}^{2\ell+1}),SU_{q}(\ell+1),\Delta)$. We denote the representation of
$C(S_{q}^{2\ell+1})$ on $L^{2}(S_{q}^{2\ell+1})$ by $\pi_{eq}$. In \cite{PsPal2}
$SU_{q}(\ell+1)$
equivariant spectral triples for this covariant representation were studied and
a non-trivial one was constructed. It is proved in \cite{SunPal} that the
Hilbert space $L^{2}(S_{q}^{2\ell+1})$ is unitarily equivalent to
 $\ell^{2}(\mathbb{N}^{\ell}\times \mathbb{Z} \times \mathbb{N}^{\ell})$
. Then the selfadjoint operator $D_{eq}$ constructed in
\cite{PsPal2} is given on the orthonormal basis $\{e_{\gamma}:\gamma \in
\mathbb{N}^{\ell}\times \mathbb{Z} \times \mathbb{N}^{\ell}\}$ by the  formula
$D_{eq}(e_{\gamma}):=d_{\gamma}e_{\gamma}$ where $d_{\gamma}$ is given by 
\begin{displaymath} 
\begin{array}{lll}
d_{\gamma}:=&\left\{\begin{array}{ll}
                            \sum_{i=1}^{2\ell +1}|\gamma_{i}| & ~\text{ if  }~~ 
(\gamma_{\ell+1},\gamma_{\ell+2},\cdots,\gamma_{2\ell+1})=0 ~~\text{and}~~
\gamma_{\ell +1}
\geq 0 \\
                            -\sum_{i=1}^{2\ell+1}|\gamma_{i}| & ~~\text{else~}\\
                       \end{array} \right.
\end{array}
\end{displaymath}
In \cite{SunPal},  a smooth subalgebra $C^{\infty}(S_{q}^{2\ell+1})
\subset C(S_{q}^{2\ell+1})$ is defined and it is shown that the spectral triple
$(C^{\infty}(S_{q}^{2\ell+1}), L^{2}(S_{q}^{2\ell+1}), D_{eq})$ is a regular
spectral triple with simple dimension spectrum $\{1,2,\cdots,2\ell+1\}$. Now we
show that the spectral triple $(C^{\infty}(S_{q}^{2\ell+1}),
L^{2}(S_{q}^{2\ell+1}), D_{eq})$ has the topological weak heat kernel
expansion. 

 We use the same notations as in \cite{SunPal}. Let $\cla_{\ell}^{\infty}:=
\Sigma^{2\ell}(C^{\infty}(\mathbb{T}))$. It follows from Corollary 4.2.3 that
$C^{\infty}(S_{q}^{2\ell+1}) \subset \cla_{\ell}^{\infty}$. let
$(C^{\infty}(S_{q}^{2\ell+1}),\pi_{\ell}, \clh_{\ell}, D_{\ell})$ be the torus
equivariant spectral triple. Let $N$ be the number operator on
$\ell^{2}(\bbn^{\ell})$ defined by 
\begin{equation*}
 Ne_{\gamma}:= (\sum_{i=1}^{\ell}\gamma_{i})e_{\gamma}
\end{equation*}

Let us denote the Hilbert space $\ell^{2}(S_{q}^{2\ell+1})$ by $\clh$. We
identify $\clh_{\ell}:=\ell^{2}(\mathbb{N}^{\ell}\times \bbz)$ with the subspace
$\ell^{2}(\bbn^{\ell}\times \bbz \times \{0\})$ and we denote the orthogonal
complement in $\clh$ by $\clh_{\ell}^{'}$. Then
$\ell^{2}(S_{q}^{2\ell+1})=\clh_{\ell}\oplus \clh_{\ell}^{'}$. Define the 
unbounded operator $D_{torus}$ on $\clh$ by the equation 
\[
D_{torus}:=  \begin{bmatrix}
                    D_{\ell} & 0 \\
                     0        & -|D_{\ell}|\otimes 1 - 1\otimes N 
                   \end{bmatrix}
\]
 Then in \cite{SunPal}, it is shown that $D_{eq}=D_{torus}$. We denote
representation $\pi_{\ell}\oplus (\pi_{e\ll}\otimes 1)$ of $C(S_{q}^{2\ell+1})$
on $\clh$ by $\pi_{torus}$. 

Let $\mathscr{T}^{\infty}:=\Sigma_{smooth}^{2}(\mathbb{C})$ and let
$\mathscr{T}_{\ell}^{\infty}:= \mathscr{T}^{\infty} \hat{\otimes}
\mathscr{T}^{\infty}
\hat{\otimes} \cdots \hat{\otimes} \mathscr{T}^{\infty}$ denote the Fr\'echet
tensor
product of $\ell$ copies.
The main theorem in \cite{SunPal} is the following.
\begin{thm}
 For every $a \in C^{\infty}(S_{q}^{2\ell+1})$, the difference
$\pi_{eq}(a)-\pi_{torus}(a) \in OP^{-\infty}_{D_{\ell}}\hat{\otimes}
\mathscr{T}_{\ell}$.
\end{thm}

Let $P_{\ell}:=\frac{1+F_{\ell}}{2}$ where $F_{\ell}:=Sign(D_{\ell})$. We denote
the rank one projection $|e_{0}\rangle \langle e_{0}|$ on
$\ell^{2}(\mathbb{N}^{\ell})$ by $P$ where $e_{0}:=e_{(0,0,\cdots,0)}$. 

\begin{ppsn}
 The equivariant spectral triple $(C^{\infty}(S_{q}^{2\ell+1}),\clh, D_{eq})$
has the topological weak heat kernel expansion.
\end{ppsn}
\textit{Proof.} Let $\mathcal{J}:=OP^{-\infty}_{D_{\ell}}\hat{\otimes}
\mathscr{T}_{\ell}^{\infty}$. In \cite{SunPal}, the following algebra is
considered.
 \begin{align*}
\mathcal{B}:=& \{a_{1}P_\ell\otimes P + a_{2}P_\ell\otimes (1-P) +
a_{3}(1-P_\ell)\otimes P +
a_{4}(1-P_\ell)\otimes (1-P)+R :\\
            &   a_{1},a_{2},a_{3},a_{4} \in \mathcal{A}_{\ell}^{\infty},  R \in
\mathcal{J}
\} \end{align*} 
The algbera $\mathcal{B}$ is isomorphic to $\cla_{\ell}^{\infty}\oplus
\cla_{\ell}^{\infty}\oplus \cla_{\ell}^{\infty}\oplus
\cla_{\ell}^{\infty}\oplus \mathcal{J}$. We give $\mathcal{B}$ the Fr\'echet
space structure coming from this decomposition. In \cite{SunPal}, it is shown
that $\mathcal{B}$ contains $C^{\infty}(S_{q}^{2\ell+1})$ and is closed under
$\delta:=[|D_{eq}|,.]$ and $d:=[D,.]$. Moreover it is shown that $\delta$ and
$d$ are continuous on $\mathcal{B}$. Note that $F_{eq}:=F_{\ell}\otimes P -
1 \otimes (1-P)$. Hence by definition $F_{eq} \in \mathcal{B}$. Now note that
the torus equivariant spectral triple $(\mathcal{A}_{\ell}^{\infty},\clh_{\ell},
D_{\ell})$ has the topological weak heat kernel asymptotic expansion. Thus it
is enough to show that the map $\tau_{2\ell+1}: (0,\infty)\times \mathcal{J}
\to \mathbb{C}$ defined by
$\tau_{2\ell+1}(t,b):=t^{2\ell+1}Tr(be^{-t|D_{eq}|})$ has uniform asymptotic
expansion. 

But this follows from the fact that $(OP^{-\infty}_{D_{\ell}},\clh_{\ell},
D_{\ell})$ and $(\mathscr{T}^{\infty}, \ell^{2}(\bbn), N)$ have the topological
weak
heat kernel expansion and by using Lemma \ref{tensor product}. This completes
the proof. \hfill $\Box$

\begin{rmrk}
 The method in \cite{SunPal} can be applied to show that the equivariant
spectral triple on the quantum $SU(2)$ constructed in \cite{PsPal_SU_q_2} has
the heat kernel asymptotic expansion property with dimension 3 and hence
deducing the dimension spectrum computed in \cite{ConSU_q2}. It has been shown
in \cite{PsPal_Dirac_SU_q_2} that the isospectral triple studied in
\cite{D-G-L-V_SU_q_2} differs from the equivariant one (with multiplicity 2)
constructed in \cite{PsPal_SU_q_2} only be a smooth perturbation. As a result it
will follow that (Since the extension $\mathcal{B}^{\infty}$ for the equivariant
spectral triple satisfying Definition \ref{topological heat kernel} contains the
algebra of
smoothing operators) the isospectral spectral triple also has the weak heat
kernel
expansion with dimension 3.
\end{rmrk}

\nocite{Higson_localindex}
\nocite{Con-Mos}
\nocite{Sch-1992}

\bibliography{references}

\def\cprime{$'$} \def\cprime{$'$}
\begin{thebibliography}{10}

\bibitem{Berline}
Nicole Berline, Ezra Getzler, and Mich{\`e}le Vergne.
\newblock {\em Heat kernels and {D}irac operators}.
\newblock Grundlehren Text Editions. Springer-Verlag, Berlin, 2004.
\newblock Corrected reprint of the 1992 original.

\bibitem{PsPal_SU_q_2}
Partha~Sarathi Chakraborty and Arupkumar Pal.
\newblock Equivariant spectral triples on the quantum {${\rm SU}(2)$} group.
\newblock {\em $K$-Theory}, 28(2):107--126, 2003.

\bibitem{PsPal_Dirac_SU_q_2}
Partha~Sarathi Chakraborty and Arupkumar Pal.
\newblock On equivariant {D}irac operators for {${\rm SU}_q(2)$}.
\newblock {\em Proc. Indian Acad. Sci. Math. Sci.}, 116(4):531--541, 2006.

\bibitem{PsPal1}
Partha~Sarathi Chakraborty and ArupKumar Pal.
\newblock Torus equivariant spectral triples for odd-dimensional quantum
  spheres coming from ${C}^{*}$ extensions.
\newblock {\em Letters in Mathematical Physics}, 80(1):57--68, April 2007.

\bibitem{PsPal2}
Partha~Sarathi Chakraborty and Arupkumar Pal.
\newblock Characterization of ${SU}_q(\ell+1)$-equivariant spectral triples for
  the odd dimensional quantum spheres.
\newblock {\em J. Reine Angew. Math.}, 623:25--42, 2008.

\bibitem{Con-Mos}
A.~Connes and H.~Moscovici.
\newblock The local index formula in noncommutative geometry.
\newblock {\em Geom. Funct. Anal.}, 5(2):174--243, 1995.

\bibitem{Con-Mos-Hopf}
A.~Connes and H.~Moscovici.
\newblock Hopf algebras, cyclic cohomology and the transverse index theorem.
\newblock {\em Comm. Math. Phys.}, 198(1):199--246, 1998.

\bibitem{ConSU_q2}
Alain Connes.
\newblock Cyclic cohomology, quantum group symmetries and the local index
  formula for {${\rm SU}\sb q(2)$}.
\newblock {\em J. Inst. Math. Jussieu}, 3(1):17--68, 2004.

\bibitem{Higson}
Nigel Higson.
\newblock Meromorphic continuation of zeta functions associated to elliptic
  operators.
\newblock In {\em Operator algebras, quantization, and noncommutative
  geometry}, volume 365 of {\em Contemp. Math.}, pages 129--142. Amer. Math.
  Soc., Providence, RI, 2004.

\bibitem{Higson_localindex}
Nigel Higson.
\newblock The residue index theorem of {C}onnes and {M}oscovici.
\newblock In {\em Surveys in noncommutative geometry}, volume~6 of {\em Clay
  Math. Proc.}, pages 71--126. Amer. Math. Soc., Providence, RI, 2006.

\bibitem{Hong-Sym1}
Jeong~Hee Hong and Wojciech Szymanski.
\newblock Quantum spheres and projective spaces as graph algebra.
\newblock {\em Commun.Math.Phys.}, 232:157--188, 2002.

\bibitem{Hong-Sym2}
Jeong~Hee Hong and Wojciech Szymanski.
\newblock Noncommutative balls and mirror quantum spheres.
\newblock {\em J.London.Math.Soc.}, March 2008.

\bibitem{SunPal}
ArupKumar Pal and S.Sundar.
\newblock Regularity and dimension spectrum of the equivariant spectral triple
  for the odd dimensional quantum spheres.
\newblock {\em Journal of Noncommutative geometry(to appear)}.

\bibitem{Sch-1992}
Larry~B. Schweitzer.
\newblock A short proof that {$M\sb n(A)$} is local if {$A$} is local and
  {F}r\'echet.
\newblock {\em Internat. J. Math.}, 3(4):581--589, 1992.

\bibitem{Sch-1993}
Larry~B. Schweitzer.
\newblock Spectral invariance of dense subalgebras of operator algebras.
\newblock {\em Internat. J. Math.}, 4(2):289--317, 1993.

\bibitem{D-G-L-V_SU_q_2}
Walter van Suijlekom, Ludwik D{\c{a}}browski, Giovanni Landi, Andrzej Sitarz,
  and Joseph~C. V{\'a}rilly.
\newblock The local index formula for {${\rm SU}_q(2)$}.
\newblock {\em $K$-Theory}, 35(3-4):375--394 (2006), 2005.

\end{thebibliography}
\bibliographystyle{plain}
\end{document}